\documentclass[reqno]{amsart}
\usepackage{amsmath,amsrefs,amssymb,tensor}

\usepackage[colorlinks=true, pdfborder={0 0 0}]{hyperref}

\numberwithin{equation}{section}

\DeclareMathOperator{\Scal}{S}
\DeclareMathOperator{\Ricci}{Ric}
\DeclareMathOperator{\Riemann}{R}
\DeclareMathOperator{\Weyl}{W}
\DeclareMathOperator{\Schouten}{P}
\DeclareMathOperator{\Bach}{B}
\DeclareMathOperator{\Beta}{B}
\DeclareMathOperator{\Sym}{Sym}
\DeclareMathOperator{\Vol}{Vol}

\DeclareMathOperator{\bigO}{O}
\DeclareMathOperator{\smallo}{o}
\newcommand{\R}{\mathbb{R}}
\renewcommand{\S}{\mathbb{S}}
\newcommand{\N}{\mathbb{N}}

\renewcommand{\[}{\left[}
\renewcommand{\]}{\right]}
\renewcommand{\(}{\left(}
\renewcommand{\)}{\right)}

\newtheorem{theorem}{Theorem}[section]
\newtheorem{proposition}{Proposition}[section]
\newtheorem{corollary}{Corollary}[section]

\newtheorem{step}{Step}[section]

\begin{document}

\title[Higher-order $Q$-curvature equation]{Existence results for the higher-order $Q$-curvature equation}

\author{Saikat Mazumdar}

\address{Saikat Mazumdar, Department of Mathematics, Indian Institute of Technology Bombay, Mumbai 400076, India}
\email{saikat@math.iitb.ac.in, saikat.mazumdar@iitb.ac.in}

\author{J\'er\^ome V\'etois}

\address{J\'er\^ome V\'etois, Department of Mathematics and Statistics, McGill University, 805 Sherbrooke Street West, Montreal, Quebec H3A 0B9, Canada}
\email{jerome.vetois@mcgill.ca}

\thanks{The second author was supported by the Discovery Grant RGPIN-2016-04195 from the Natural Sciences and Engineering Research Council of Canada. This work was initiated when the first author held a postdoctoral position at McGill University under the co-supervision of Professors Pengfei Guan, Niky Kamran and the second author, that was supported by the NSERC Discovery Grants RGPIN-04443-2018, RGPIN-05490-2018 and RGPIN-04195-2016.}

\date{December 21, 2022}

\begin{abstract}
We obtain existence results for the $Q$-curvature equation of order $2k$ on a closed Riemannian manifold of dimension $n\ge 2k+1$, where $k\ge1$ is an integer. We obtain these results under the assumptions that the Yamabe invariant of order $2k$ is positive and the Green's function of the corresponding operator is positive, which are satisfied in particular when the manifold is Einstein with positive scalar curvature. In the case where $2k+1\le n\le2k+3$ or the manifold is locally conformally flat, we assume moreover that the operator has positive mass. In the case where $n\ge2k+4$ and the manifold is not locally conformally flat, the results essentially reduce  to the determination of the sign of a complicated constant depending only on $n$ and $k$.
\end{abstract}

\maketitle

\section{Introduction and main results}\label{Sec1}

Given an integer $k\ge1$, a smooth, closed Riemannian manifold $\(M,g\)$ of dimension $n>2k$ and a smooth positive function $f$ in $M$, we consider the equation
\begin{equation}\label{Eq1}
P_{2k}u=f\left|u\right|^{2^*_k-2}u\quad\text{in M},
\end{equation}
where $P_{2k}$ is the GJMS operator with leading part $\Delta^k$, $\Delta:=\delta d$ is the Laplace--Beltrami operator with nonnegative eigenvalues and $2^*_k:=2n/\(n-2k\)$ is the critical Sobolev exponent. The so-called GJMS operators were discovered by Graham, Jenne, Mason and Sparling~\cite{GraJenMasSpa} by using a construction based on the Fefferman--Graham ambient metric~\cites{FefGra1,FefGra2}. They provide a natural extension to higher orders of the Yamabe operator~\cite{Yam} ($k=1$) and the Paneitz--Branson operator~\cites{Bra1,Pan} ($k=2$). When $u$ is positive, \eqref{Eq1} arises in the problem of prescribing \nobreak Branson's $Q$-curvature of order $2k$ in a given conformal class (see Branson~\cites{Bra2}). More precisely, the positive solutions $u$ to the equation \eqref{Eq1} correspond to the conformal metrics $u^{4/\(n-2k\)}g$ with $Q$-curvature of order $2k$ equal to $\frac{2}{n-2k}f$. 

\smallskip
Let $Y_{2k}$ be the conformal invariant defined by
$$Y_{2k}:=\inf_{\widetilde{g}\in\[g\]}\(\Vol_{\widetilde{g}}\(M\)^{-\frac{n-2k}{n}}\int_MQ_{2k,\widetilde{g}}\,dv_{\widetilde{g}}\)=\inf_{\substack{u\in C^\infty\(M\)\\u>0\text{ in }M}}\frac{\displaystyle\int_MuP_{2k}u\,dv_g}{\displaystyle\(\int_Mu^{2^*_k}dv_g\)^{\frac{n-2k}{n}}},$$
where $\[g\]$ is the conformal class of $g$, and $\Vol_{\widetilde{g}}\(M\)$, $dv_{\widetilde{g}}$ and $Q_{2k,\widetilde{g}}$ are the volume, volume element and $Q$-curvature of order $2k$, respectively, of $\(M,\widetilde{g}\)$. Throughout this paper, we assume that $Y_{2k}>0$. As is easily seen, this is equivalent to the coercivity of the operator $P_{2k}$, which is also equivalent to $\lambda_1\(P_{2k}\)>0$, where $\lambda_1\(P_{2k}\)$ is the first eigenvalue of $P_{2k}$.

\smallskip
In the case where $k=1$, it is well-known that there exists at least one positive solution to the equation \eqref{Eq1} with $f\equiv1$ if and only if $Y_2>0$ (see the historic work of Aubin~\cite{Aub}, Schoen~\cite{Sch1}, Trudinger~\cite{Tru} and Yamabe~\cite{Yam}). In the case where $k=2$, the existence of at least one positive solution to this problem has been obtained under positivity assumptions on the scalar curvature and $Q$-curvature of order 4 (see Gursky and Malchiodi~\cite{GurMal}) and later extended to the cases where $Y_2>0$ and $Y_4>0$ in dimension $n\ge6$ (see Gursky, Hang and Lin~\cites{GurHangLin}) and the case where $Y_2>0$ and $Q_4>0$ in dimension $n\ge5$ (see Hang and Yang~\cites{HangYang1,HangYang2}). This question has also been solved by Qing and Raske~\cite{QingRas} in the locally conformally flat case for all orders $k\ge2$, under a topological assumption on the Poincar\'e exponent of the holonomy representation of the fundamental group, using an approach introduced by Schoen~\cite{Sch2} for $k=1$. More general existence results have also been obtained in the case where $f\not\equiv1$ (see among others Aubin~\cite{Aub}, Escobar and Schoen~\cite{EscSch}, Hebey~\cite{Heb} and Hebey and Vaugon~\cite{HebVau} for $k=1$, Djadli, Hebey and Ledoux~\cite{DjaHebLed}, Esposito and Robert~\cite{EspRob} and Robert~\cite{Rob1} for $k=2$, Chen and Hou~\cite{ChenHou} for $k=3$ and Robert~\cite{Rob2} for higher orders). 

\smallskip
We let $\Weyl$ be the Weyl tensor of $\(M,g\)$ and $\left|\Weyl\right|$ be the norm of $\Weyl$ with respect to $g$. In the case where $2k+1\le n\le 2k+3$ or $\(M,g\)$ is locally conformally flat, assuming that $Y_{2k}>0$, for every point $\xi\in M$, we let $m\(\xi\)$ be the mass of $P_{2k}$ at $\xi$ (see \eqref{Sec3Eq2} for the definition of the mass). Our main result is the following:

\begin{theorem}\label{Th1}
Let $k\ge1$ be an integer, $\(M,g\)$ be a smooth, closed Riemannian manifold of dimension $n\ge 2k+1$ and $f$ be a smooth positive function in $M$. Assume that $Y_{2k}>0$ and there exists a maximal point $\xi$ of $f$ such that 
\begin{equation}\label{Th1Eq1}
\Delta f\(\xi\)=0\quad\text{if }n\ge 2k+2
\end{equation}
and 
\begin{equation}\label{Th1Eq2}
\left\{\begin{aligned}&\left|\Weyl\(\xi\)\right|^2f\(\xi\)+c\(n,k\)\Delta^2f\(\xi\)>0&&\text{if }n\ge2k+5\\
&\Weyl\(\xi\)\ne0&&\text{if }n=2k+4\\
&m\(\xi\)>0&&\text{if }2k+1\le n\le 2k+3,
\end{aligned}\right.
\end{equation}
where $c\(n,k\)$ is a positive constant depending only on $n$ and $k$ (see \eqref{Th1Eq5} for the value of $c\(n,k\)$). Then there exists a nontrivial solution $u\in C^{2k}\(M\)$ to the equation \eqref{Eq1}, which minimizes the energy functional \eqref{Eq2}. If moreover the Green's function of the operator $P_{2k}$ is positive, then $u$ is positive, which implies that the $Q$-curvature of order $2k$ of the metric $u^{4/\(n-2k\)}g$ is equal to $\frac{2}{n-2k}f$. 
\end{theorem}

In particular, Theorem~\ref{Th1} extends to all orders previous results obtained by Aubin~\cite{Aub} for $k=1$ (in this case, the positivity of the Green's function is not an issue), Esposito and Robert~\cite{EspRob} for $k=2$ and Chen and Hou~\cite{ChenHou} for $k=3$.

\smallskip
In the case where $f$ is constant, we obtain the following:

\begin{theorem}\label{Th2}
Let $k\ge1$ be an integer and $\(M,g\)$ be a smooth, closed Riemannian manifold of dimension $n\ge 2k+1$. Assume that $Y_{2k}>0$ and its Green's function is positive. Assume moreover that if $2k+1\le n\le 2k+3$ or $\(M,g\)$ is locally conformally flat, then $m\(\xi\)>0$ for some point $\xi\in M$. Then there exists a conformal metric to $g$ with constant $Q$-curvature of order $2k$. 
\end{theorem}

Notice that Theorem~\ref{Th2} is a direct consequence of Theorem~\ref{Th1} in the case where $\(M,g\)$ is not locally conformally flat of dimension $n\ge2k+4$. A more general result about the locally conformally flat case will be stated in Section~\ref{Sec3}.

\smallskip
When $\(M,g\)$ is Einstein, Fefferman and Graham~\cite{FefGra2}*{Proposition~7.9} (see also Gover~\cite{Gov} for a proof based on tractors) established the formula
$$P_{2k}=\prod_{j=1}^k\(\Delta+\frac{\(n+2j-2\)\(n-2j\)}{4n\(n-1\)}\Scal\),$$
where $\Scal$ is the Scalar curvature of $\(M,g\)$. In this case, it is easy to see that if $\Scal$ is positive, then $P_{2k}$ is coercive, and so $Y_{2k}>0$. Furthermore, successive applications of the maximum principles yield that the Green's function of the operator $P_{2k}$ is positive. Therefore, we obtain the following corollary of Theorem~\ref{Th1}:

\begin{corollary}
Let $k\ge1$ be an integer and $\(M,g\)$ be a smooth, closed Einstein manifold of positive scalar curvature and dimension $n\ge 2k+1$. Let $f$ be a smooth positive function in $M$ such that there exists a maximal point $\xi$ of $f$ satisfying \eqref{Th1Eq1} and \eqref{Th1Eq2}. Then there exists a conformal metric to $g$ with $Q$-curvature of order $2k$ equal to~$\frac{2}{n-2k}f$.
\end{corollary}

The positivity of the Green's function of the operator $P_4$ has been shown to be true by  Gursky and Malchiodi~\cites{GurMal} and Hang and Yang~\cites{HangYang1,HangYang2} under positivity assumptions on the $Q$-curvature of order 4 and the scalar curvature or the Yamabe invariant of the manifold. Positivity results for the mass of $P_4$ have also been obtained by Gursky and Malchiodi~\cites{GurMal}, Hang and Yang~\cite{HangYang1}, Humbert and Raulot~\cite{HumRau} and Michel~\cite{Mic}, thus extending the positive mass theorem obtained by Schoen and Yau~\cites{SchYau1,SchYau2,SchYau3} for $k=1$. As far as the authors know, no such results have yet been obtained for higher orders. As regards the case where $n=2k$, we point out that the problem of prescribing the $Q$-curvature involves a different equation than \eqref{Eq1} which contains an exponential non-linearity. Some references in this case are Chang and Yang~\cite{ChangYang}, Djadli and Malchiodi~\cite{DjaMal} and Li, Li and Liu~\cite{LiLiLiu} for $k=2$ and Baird, Fardoun and Regbaoui~\cite{BaiFarReg} for higher orders.

\smallskip
The proofs of Theorems~\ref{Th1} and~\ref{Th2} are based on the approach introduced by Aubin~\cite{Aub} and Schoen~\cite{Sch1} in the case where $k=1$. This approach consists in deriving an asymptotic expansion for the energy functional associated with the equation \eqref{Eq1}, which we apply to a suitable family of test functions depending on a real parameter (see \eqref{Eq2} for the energy functional; see \eqref{Eq6} and \eqref{Sec3Eq4} for the definitions of our families of test functions). To simplify the calculations of curvature terms, we use the conformal normal coordinates introduced by Lee and Parker~\cite{LeePar} and later improved by Cao~\cite{Cao} and G\"{u}nther~\cite{Gun}. Our proof also crucially relies on the derivation of an expression for the highest-order terms of the GJMS operators (see \eqref{Pr1Step1Eq1}), which we obtain by using Juhl's formulae~\cite{Juhl}. In the case where $n\ge2k+4$, the proof essentially reduces to determining the sign of a constant $C\(n,k\)$, which appears in the energy expansion (see \eqref{Pr1Eq1}). In particular, we recover the values found in~\cites{ChenHou,EspRob} for $C\(n,k\)$ with $k\in\left\{2,3\right\}$. We then conclude the proof by using a minimization result in the spirit of Aubin~\cite{Aub} (see Mazumdar~\cite{Maz}*{Theorem~3}). When the Green's function of the operator $P_{2k}$ is positive, by an application of the Green's representation formula, we obtain moreover that the minimizing solution is positive (see the argument in~\cite{Maz}*{end of Section 3}). We point out that at one place in the proof, namely in the very last computation to determine the sign of $C\(n,k\)$ (see \eqref{Pr1Eq11}), we have used the computation software {\it Maple} to expand a complicated polynomial with integer coefficients.

\smallskip
The paper is organized as follows. In Section~\ref{Sec2}, we prove Theorem~\ref{Th1} in the case where $n\ge 2k+4$. In Section~\ref{Sec3}, we complete the proof of Theorems~\ref{Th1} in the remaining case where $2k+1\le n\le 2k+3$ and we state and prove a more general result in the case where $g$ is conformally flat in some open subset of the manifold. Theorem~\ref{Th2} then directly follows from this new result together with Theorem~\ref{Th1}. 

\section{Proof of Theorem~\ref{Th1} in the case where $n\ge 2k+4$}\label{Sec2}

Given an integer $k\ge1$ and a smooth positive function $f$ in $M$, we let $I_{k,f}$ be the energy functional defined as
\begin{equation}\label{Eq2}
I_{k,f,g}\(u\):=\frac{\displaystyle\int_MuP_{2k}u\,dv_g}{\displaystyle\(\int_Mf\left|u\right|^{2^*_k}dv_g\)^{\frac{n-2k}{n}}}
\end{equation}
for all functions $u\in C^{2k}\(M\)$ such that $u\not\equiv0$. We fix a point $\xi\in M$. By applying a conformal change of metric (see Cao~\cite{Cao}, G\"{u}nther~\cite{Gun} and Lee and Parker~\cite{LeePar}), we may assume that 
\begin{equation}\label{Eq3}
\det g\(x\)=1\quad\forall x\in\Omega
\end{equation}
for some neighborhood $\Omega$ of the point $\xi$, where $\det g$ is the determinant of $g$ in geodesic normal coordinates at $\xi$. In particular (see~\cite{LeePar}), it follows from \eqref{Eq3} that
\begin{equation}\label{Eq4}
\Ricci\(\xi\)=\Sym\nabla\Ricci\(\xi\)=\Sym\(\Ricci_{ab;cd}\(\xi\)+\frac{2}{9}\Weyl_{eabf}\(\xi\)\tensor{\Weyl}{^e_{cd}^f}\(\xi\)\)=0,
\end{equation}
where $\Sym$ stands for the symmetric part, $\Ricci$ is the Ricci tensor, and $\Ricci_{ab;cd}$ and $\Weyl_{eabf}$ are the coordinates of $\nabla^2\Ricci$ and $\Weyl$, respectively, with the standard convention on raising and lowering indices. By taking traces in \eqref{Eq4} and using Bianchi's identities, we obtain
\begin{equation}\label{Eq5}
\Scal\(\xi\)=\left|\nabla \Scal\(\xi\)\right|=0,\ \Delta \Scal\(\xi\)=\frac{1}{6}\left|\Weyl\(\xi\)\right|^2\text{ and }\tensor{\Ricci}{^{ab}_{;ab}}\(\xi\)=-\frac{1}{12}\left|\Weyl\(\xi\)\right|^2
\end{equation}
Let $r_0>0$ be such that the injectivity radius of the metric $g$ at the point $\xi$ is greater than $3r_0$ and $B\(\xi,3r_0\)\subset\Omega$, where $B\(\xi,r_0\)$ is the ball of center $\xi$ and radius $3r_0$ with respect to $g$. We then let $\chi$ be a smooth cutoff function in $\[0,\infty\)$ such that $\chi\equiv1$ in $\[0,r_0\]$, $0\le\chi\le1$ in $\(r_0,2r_0\)$ and $\chi\equiv0$ in $\[2r_0,\infty\)$. For every $\mu>0$, we then define our test functions as
\begin{equation}\label{Eq6}
U_\mu\(x\):=\chi\(d_g\(x,\xi\)\)\mu^{\frac{2k-n}{2}}U\big(\mu^{-1}\exp_\xi^{-1}x\big)\quad\forall x\in M,
\end{equation}
where $d_g$ is the geodesic distance with respect to $g$, $\exp_\xi$ is the exponential map with respect to $g$ at the point $\xi$ and $U$ is the function in $\R^n$ (we identify $T_\xi M$ with $\R^n$) defined as
$$U\(x\):=\big(1+\left|x\right|^2\big)^{-\frac{n-2k}{2}}\quad\forall x\in\R^n.$$
It is easy to verify that $U$ is a solution of the equation 
$$\Delta_0^k\,U=\[\prod_{j=-k}^{k-1}\(n+2j\)\]U^{2^*_k-1}\quad\text{in }\R^n,$$
where $\Delta_0$ is the Euclidean Laplacian.

\begin{proposition}\label{Pr1}
Let $k\ge1$ be an integer, $\(M,g\)$ be a smooth, closed Riemannian manifold of dimension $n\ge 2k+4$ and $f$ be a smooth positive function in $M$. Assume that $g$ satisfies \eqref{Eq3} for some point $\xi\in M$. Let $I_{k,f,g}$ be as in \eqref{Eq2} and $U_\mu$ be as in \eqref{Eq6}. Then there exists a positive constant $C\(n,k\)$ depending only on $n$ and $k$ (see \eqref{Pr1Eq10} for the value of $C\(n,k\)$) such that as $\mu\to0$,
\begin{align}\label{Pr1Eq1}
I_{k,f,g}\(U_\mu\)&=\omega_n^{\frac{2k}{n}}f\(\xi\)^{-\frac{n-2k}{n}}\Bigg(\(2k-1\)!\Beta\(\frac{n}{2}-k,2k\)^{-1}\nonumber\\
&\quad\times\Bigg(1+\frac{\(n-2k\)}{2n\(n-2\)}\frac{\Delta f\(\xi\)}{f\(\xi\)}\mu^2\nonumber\\
&\quad-\frac{\(n-2k\)}{4n\(n-2\)}\(\frac{\Delta^2f\(\xi\)}{2\(n-4\)f\(\xi\)}-\frac{\(n-k\)\(\Delta f\(\xi\)\)^2}{n\(n-2\)f\(\xi\)^2}\)\mu^4\Bigg)\nonumber\\
&\quad-C\(n,k\)\,\mu^4\left\{\begin{aligned}&\left|\Weyl\(\xi\)\right|^2\ln\(1/\mu\)+\bigO\(1\)&&\text{if }n=2k+4\\&\left|\Weyl\(\xi\)\right|^2+\smallo\(1\)&&\text{if }n>2k+4.\end{aligned}\right\}\Bigg),
\end{align}
where $\omega_n$ is the volume of the standard $n$-dimensional sphere and $\Beta$ is the beta function defined as 
$$\Beta\(a,b\)=\frac{\Gamma\(a\)\Gamma\(b\)}{\Gamma\(a+b\)}\quad\forall a,b>0.$$
\end{proposition}

\proof[Proof of Proposition~\ref{Pr1}]
We let $\Schouten$ be the Schouten tensor defined as
$$\Schouten:=\frac{1}{n-2}\(\Ricci-\frac{\Scal}{2\(n-1\)}\,g\)$$
and $\Bach$ be the Bach tensor whose coordinates are given by
$$\Bach_{ij}:=\Schouten_{ab}\tensor{\Weyl}{_i^a_j^b}+\tensor{\Schouten}{_{ij;a}^a}-\tensor{\Schouten}{_{ia;j}^a},$$
where $\Weyl_{iajb}$, $\Schouten_{ab}$ and $\Schouten_{ij;ab}$ are the coordinates of $\Weyl$, $\Schouten$ and $\nabla^2\Schouten$, respectively. We let $\(\cdot,\cdot\)$ be the multiple inner product induced by the metric $g$ for the tensors of same rank, i.e. such that $\(S,T\)=S^{i_1\dotsc i_l}T_{i_1\dotsc i_l}$ for all tensors $S$ and $T$ of rank $l\in\N$. The first step in the proof of Proposition~\ref{Pr1} is as follows:

\begin{step}\label{Pr1Step1}
For every $k\in\N$ such that $n\ge 2k+1$, we have
\begin{align}\label{Pr1Step1Eq1}
P_{2k}&=\Delta^k+k\Delta^{k-1}\(J_1\cdot\)+k\(k-1\)\Delta^{k-2}\(J_2\cdot+\(T_1,\nabla\)+\(T_2,\nabla^2\)\)\nonumber\\
&\quad+k\(k-1\)\(k-2\)\Delta^{k-3}\(\(T_3,\nabla^2\)+\(T_4,\nabla^3\)\)\allowdisplaybreaks\nonumber\\
&\quad+k\(k-1\)\(k-2\)\(k-3\)\Delta^{k-4}\(T_5,\nabla^4\)+Z,
\end{align}
where $Z$ is a smooth linear operator of order less than $2k-4$ if $k\ge3$, $Z:=0$ if $k\le2$, the functions $J_1$ and $J_2$ are defined as 
$$J_1:=\frac{n-2}{4\(n-1\)}\Scal$$
and
$$J_2:=\frac{1}{6}\(\frac{3n^2-12n-4k+8}{16\(n-1\)^2}\Scal^2-\(k+1\)\(n-4\)\left|\Schouten\right|^2-\frac{3n+2k-4}{4\(n-1\)}\Delta\Scal\),$$
and the tensors $T_1$, $T_2$, $T_3$, $T_4$ and $T_5$ are defined as
\begin{align*}
&T_1:=\frac{n-2}{4\(n-1\)}\nabla\Scal-\frac{2}{3}\(k+1\)\delta \Schouten,\allowdisplaybreaks\\
&T_2:=\frac{2}{3}\(k+1\)\Schouten,\allowdisplaybreaks\\
&T_3:=\frac{n-2}{6\(n-1\)}\nabla^2\Scal+\frac{\(k+1\)\(n-2\)}{6\(n-1\)}\Scal\Schouten-\frac{k+1}{3}\(\delta\nabla\Schouten+2\nabla\delta\Schouten+2\Riemann\ast\Schouten\)\\
&\qquad-\frac{2}{15}\(k+1\)\(k+2\)\(3\Schouten^{\#}\Schouten+\frac{\Bach}{n-4}\),\allowdisplaybreaks\\
&T_4:=\frac{2}{3}\(k+1\)\nabla \Schouten
\end{align*}
and
$$T_5:=\frac{2}{5}\(k+1\)\(\frac{5k+7}{9}\Schouten\otimes \Schouten+\nabla^2\Schouten\),$$
where $\#$ stands for the musical isomorphism with respect to $g$ (i.e. $\Schouten^\#:=g^{-1}\Schouten$), and $\delta\nabla\Schouten$, $\nabla\delta\Schouten$ and $\Riemann\ast\Schouten$ stand for the covariant tensors whose coordinates are given by 
\begin{equation}\label{Pr1Step1Eq2}
\(\delta\nabla\Schouten\)_{ij}:=-\tensor{\Schouten}{_{ij;a}^a},\ \(\nabla\delta\Schouten\)_{ij}:=-\tensor{\Schouten}{_i^a_{;aj}}\text{ and }\(\Riemann\ast\Schouten\)_{ij}:=\tensor{\Riemann}{_{ia}^a_b}\tensor{\Schouten}{_j^b}+\Riemann_{ibja}\Schouten^{ab},
\end{equation}
where $\Riemann_{ibja}$, $\Schouten_{ab}$ and $\Schouten_{ij;ab}$ are the coordinates of the Riemann tensor, $\Schouten$ and $\nabla^2\Schouten$, respectively.
\end{step}

\proof[Proof of Step~\ref{Pr1Step1}]
Throughout this proof, for every integer $l$, $\smallo^l$ stands for a linear operator of order less than $l$ if $l>0$ and $\smallo^{l}:=0$ if $l\le0$. Juhl's formulae~\cite{Juhl} (see also Fefferman and Graham~\cite{FefGra3}) give
\begin{align}\label{Pr1Step1Eq3}
&P_{2k}=M_2^k-\sum_{j=1}^{k-1}j\(k-j\)M_2^{j-1}M_4M_2^{k-j-1}\nonumber\\
&+\frac{1}{4}\sum_{j=1}^{k-2}j\(j+1\)\(k-j\)\(k-j-1\)M_2^{j-1}M_6M_2^{k-j-2}\nonumber\allowdisplaybreaks\\
&+\sum_{j=2}^{k-2}\(j+1\)\(k-j-1\)\sum_{i=1}^{j-1}i\(k-i\)M_2^{i-1}M_4M_2^{j-i-1}M_4M_2^{k-j-2}+\smallo^{2k-5},
\end{align}
where the operators $M_2$, $M_4$ and $M_6$ are defined as
$$M_2:=\Delta+\mu_2,\quad M_4:=4\delta\Schouten^\#d+\mu_4\quad\text{and}\quad M_6:=\delta A_6^\#d+\mu_6,$$
where $\mu_6$ is a smooth function in $M$ which we do not need explicitly, $\mu_2$ and $\mu_4$ are the functions defined as
$$\mu_2:=\frac{n-2}{4\(n-1\)}\Scal\quad\text{ and }\quad \mu_4:=\frac{\Delta\Scal}{2\(n-1\)}+\frac{\Scal^2}{4\(n-1\)^2}+\(n-4\)\left|\Schouten\right|^2$$
and $A_6$ is the tensor defined as
$$A_6:=48\Schouten^{\#}\Schouten+\frac{16}{n-4}\Bach.$$
We point out that throughout this paper, we use the same sign convention for the Riemann tensor as in the paper of Lee and Parker~\cite{LeePar}, which is the opposite of the convention used by Juhl~\cite{Juhl}. Straightforward expansions yield
\begin{align}\label{Pr1Step1Eq6}
&M_2^k=\Delta^k+\frac{n-2}{4\(n-1\)}\sum_{j=1}^k\Delta^{j-1}\(\Scal\Delta^{k-j}\)\nonumber\\
&\quad+\frac{\(n-2\)^2}{16\(n-1\)^2}\sum_{j=2}^k\sum_{i=1}^{j-1}\Delta^{i-1}\(\Scal\Delta^{j-i-1}\(\Scal\Delta^{k-j}\)\)+\smallo^{2k-5}\allowdisplaybreaks\nonumber\\
&=\Delta^k+\frac{n-2}{4\(n-1\)}\sum_{j=1}^k\Delta^{j-1}\(\Scal\Delta^{k-j}\)+\frac{\(n-2\)^2}{16\(n-1\)^2}\sum_{j=2}^k\(j-1\)\Delta^{k-2}\(\Scal^2\cdot\)+\smallo^{2k-4}\allowdisplaybreaks\nonumber\\
&=\Delta^k+\frac{n-2}{4\(n-1\)}\sum_{j=1}^k\Delta^{j-1}\(\Scal\Delta^{k-j}\)+\frac{k\(k-1\)\(n-2\)^2}{32\(n-1\)^2}\Delta^{k-2}\(\Scal^2\cdot\)+\smallo^{2k-4}
\end{align}
and
\begin{align}\label{Pr1Step1Eq7}
&M_2^{j-1}M_4M_2^{k-j-1}=4\Delta^{j-1}\delta\Schouten^\#d\Delta^{k-j-1}+\Delta^{j-1}\(\mu_4\Delta^{k-j-1}\)\nonumber\\
&\quad+\frac{n-2}{n-1}\sum_{i=1}^{j-1}\Delta^{i-1}\big(\Scal\Delta^{j-i-1}\delta \Schouten^\#d\Delta^{k-j-1}\big)\nonumber\\
&\quad+\frac{n-2}{n-1}\sum_{i=j+1}^{k-1}\Delta^{j-1}\delta \Schouten^\#d\Delta^{i-j-1}\(\Scal\Delta^{k-i-1}\)+\smallo^{2k-5}\allowdisplaybreaks\nonumber\\
&=4\Delta^{j-1}\delta \Schouten^\#d\Delta^{k-j-1}+\mu_4\Delta^{k-2}-\frac{\(k-2\)\(n-2\)}{n-1}\Delta^{k-3}\(\Scal\Schouten,\nabla^2\)+\smallo^{2k-4}
\end{align}
and
\begin{equation}\label{Pr1Step1Eq8}
M_2^{j-1}M_6M_2^{k-j-2}=\Delta^{j-1}\delta A_6^\#d\Delta^{k-j-2}+\smallo^{2k-5}=-\Delta^{k-3}\(A_6,\nabla^2\)+\smallo^{2k-4}
\end{equation}
and
\begin{align}\label{Pr1Step1Eq9}
M_2^{i-1}M_4M_2^{j-i-1}M_4M_2^{k-j-2}&=16\Delta^{i-1}\delta \Schouten^\#d\Delta^{j-i-1}\delta \Schouten^\#d\Delta^{k-j-2}+\smallo^{2k-5}\nonumber\\
&=16\Delta^{k-4}\(\Schouten\otimes\Schouten,\nabla^4\)+\smallo^{2k-4}.
\end{align}
Furthermore, by induction, one can check that
\begin{equation}\label{Pr1Step1Eq10}
\Scal\Delta^{j}=\Delta^j\(\Scal\cdot\)-j\(\Delta\Scal\)\Delta^{j-1}+2j\Delta^{j-1}\(\nabla\Scal,\nabla\)+2j\(j-1\)\Delta^{j-2}\(\nabla^2\Scal,\nabla^2\)+\smallo^{2j-2}
\end{equation}
and
\begin{multline}\label{Pr1Step1Eq11}
\delta \Schouten^\#d\Delta^j=\Delta^j\(\(\delta \Schouten,\nabla\)-\(\Schouten,\nabla^2\)\)+j\Delta^{j-1}\big(\(\delta\nabla\Schouten+2\nabla\delta\Schouten+2\Riemann\ast \Schouten,\nabla^2\)\\
-2\(\nabla \Schouten,\nabla^3\)\big)-2j\(j-1\)\Delta^{j-2}\(\nabla^2\Schouten,\nabla^4\)+\smallo^{2j},
\end{multline}
where $\delta\nabla\Schouten$, $\nabla\delta\Schouten$ and $\Riemann\ast \Schouten$ are as in \eqref{Pr1Step1Eq2}. The proof of \eqref{Pr1Step1Eq11} relies on the commutation formula
$$u_{;abcd}=u_{;cdab}+\tensor{\Riemann}{^e_{cad}}u_{;eb}+\tensor{\Riemann}{^e_{abd}}u_{;ce}+\tensor{\Riemann}{^e_{cbd}}u_{;ae}+\tensor{\Riemann}{^e_{abc}}u_{;de}+\smallo^2u,$$
which gives 
\begin{align*}
&\delta \Schouten^{\#}d\Delta u-\Delta\delta \Schouten^{\#}du=(\Schouten^{bc}\tensor{u}{_{;a}^a_b})_{;c}-\tensor{(\Schouten^{bc}u_{;b})}{_{;ca}^a}\\
&\ \ =\Schouten^{bc}\(\tensor{u}{_{;a}^a_{bc}}-\tensor{u}{_{;bca}^a}\)-\tensor{\Schouten}{^{bc}_{;a}^a}u_{;bc}-2\tensor{\Schouten}{^{bc}_{;ca}}\tensor{u}{_{;b}^a}-2\tensor{\Schouten}{^{bc}_{;a}}\tensor{u}{_{;bc}^a}+\smallo^2u\\
&\ \ =2\Schouten^{bc}\big(\tensor{\Riemann}{^d_b^a_c}u_{;ad}+\tensor{\Riemann}{^d_a^a_c}u_{;bd}\big)-\tensor{\Schouten}{^{bc}_{;a}^a}u_{;bc}-2\tensor{\Schouten}{^{bc}_{;ca}}\tensor{u}{_{;b}^a}-2\tensor{\Schouten}{^{bc}_{;a}}\tensor{u}{_{;bc}^a}+\smallo^2u\\
&\ \ =\(\delta\nabla\Schouten+2\nabla\delta\Schouten+2\Riemann\ast \Schouten,\nabla^2u\)-2\(\nabla \Schouten,\nabla^3u\)+\smallo^2u.
\end{align*}
By combining Faulhaber's formulae with \eqref{Pr1Step1Eq7}--\eqref{Pr1Step1Eq11}, we obtain
\begin{align}\label{Pr1Step1Eq12}
&\sum_{j=1}^k\Delta^{j-1}\(\Scal\Delta^{k-j}\)=k\Delta^{k-1}\(\Scal \cdot\)-\frac{k\(k-1\)}{2}\Delta^{k-2}\(\(\Delta\Scal\)\cdot\)\nonumber\\
&\quad+k\(k-1\)\Delta^{k-2}\(\nabla\Scal,\nabla\)+\frac{2k\(k-1\)\(k-2\)}{3                                                                                                                                                                                                                                                          }\Delta^{k-3}\(\nabla^2\Scal,\nabla^2\)+\smallo^{2k-4}
\end{align}
and
\begin{align}
&\sum_{j=1}^{k-1}j\(k-j\)M_2^{j-1}M_4M_2^{k-j-1}=k\(k-1\)\(k+1\)\bigg(\frac{2}{3}\Delta^{k-2}\(\(\delta \Schouten,\nabla\)-\(\Schouten,\nabla^2\)\)\nonumber\\
&\quad+\frac{k-2}{3}\Delta^{k-3}\(\(\delta\nabla\Schouten+2\nabla\delta\Schouten+2\Riemann\ast \Schouten,\nabla^2\)-2\(\nabla \Schouten,\nabla^3\)\)+\frac{1}{6}\Delta^{k-2}\(\mu_4\cdot\)\nonumber\\
&\quad-\frac{2\(k-2\)\(k-3\)}{5}\Delta^{k-4}\(\nabla^2\Schouten,\nabla^4\)-\frac{\(k-2\)\(n-2\)}{6\(n-1\)}\Delta^{k-3}\(\Scal\Schouten,\nabla^2\)\bigg)+\smallo^{2k-4}\label{Pr1Step1Eq13}
\end{align}
and
\begin{align}\label{Pr1Step1Eq14}
&\sum_{j=1}^{k-2}j\(j+1\)\(k-j\)\(k-j-1\)M_2^{j-1}M_6M_2^{k-j-2}\nonumber\\
&\quad=-\frac{k\(k-1\)\(k-2\)\(k+1\)\(k+2\)}{30}\Delta^{k-3}\(A_6,\nabla^2\)+\smallo^{2k-4}
\end{align}
and
\begin{align}\label{Pr1Step1Eq15}
&\sum_{j=2}^{k-2}\(j+1\)\(k-j-1\)\sum_{i=1}^{j-1}i\(k-i\)M_2^{i-1}M_4M_2^{j-i-1}M_4M_2^{k-j-2}\nonumber\\
&\quad=\frac{2k\(k-1\)\(k-2\)\(k-3\)\(k+1\)\(5k+7\)}{45}\Delta^{k-4}\(\Schouten\otimes \Schouten,\nabla^4\)+\smallo^{2k-4}.
\end{align}
Finally, \eqref{Pr1Step1Eq1} follows by putting together \eqref{Pr1Step1Eq3}, \eqref{Pr1Step1Eq6} and \eqref{Pr1Step1Eq12}--\eqref{Pr1Step1Eq15}. This ends the proof of Step~\ref{Pr1Step1}.
\endproof

The next step is as follows:

\begin{step}\label{Pr1Step2}
Assume that $n\ge 2k+4$ and $k\ge3$. Then for every smooth linear operator $Z$ of order less than $2k-4$, as $\mu\to0$,
\begin{equation}\label{Pr1Step2Eq1}
\int_MU_\mu ZU_\mu\,dv_g=\left\{\begin{aligned}&\bigO\(\mu^4\)&&\text{if }n=2k+4\\&\smallo\(\mu^4\)&&\text{if }n>2k+4.\end{aligned}\right.
\end{equation}
\end{step}

\proof[Proof of Step~\ref{Pr1Step2}]
By rewriting the integral in geodesic normal coordinates, we obtain
\begin{equation}\label{Pr1Step2Eq2}
\int_MU_\mu ZU_\mu\,dv_g=\int_{B\(0,2r_0\)}\widetilde{U}_\mu\widetilde{Z}\widetilde{U}_\mu\,dx=\sum_{\left|\alpha\right|<2k-4}\int_{B\(0,2r_0\)} z_\alpha\widetilde{U}_\mu\partial^{\alpha}\widetilde{U}_\mu\,dx,
\end{equation}
where 
\begin{equation}\label{Pr1Step2Eq3}
\widetilde{U}_\mu\(x\):=\mu^{\frac{2k-n}{2}}U\(x/\mu\)\text{ and }\widetilde{Z}\(x\):=\sum_{\left|\alpha\right|<2k-4}z_\alpha\(x\)\partial^{\(\alpha\)}\quad\forall x\in B\(0,2r_0\)
\end{equation}
for some smooth functions $z_\alpha$ in $B\(0,2r_0\)$, where $\alpha$ is a multi-index. A straightforward change of variable then gives
\begin{equation}\label{Pr1Step2Eq4}
\int_{B\(0,2r_0\)} z_\alpha\widetilde{U}_\mu\partial^{\(\alpha\)}\widetilde{U}_\mu\,dx=\mu^{2k-\left|\alpha\right|}\int_{B\(0,2r_0/\mu\)} z_\alpha\(\mu x\)U\(x\)\partial^{\(\alpha\)}U\(x\)dx.
\end{equation}
An easy induction yields that for every multi-index $\alpha$, there exists a constant $C_\alpha$ such that 
\begin{equation}\label{Pr1Step2Eq5}
\big|\partial^{\(\alpha\)}U\(x\)\big|\le C_\alpha\big(1+\left|x\right|^2\big)^{-\frac{n-2k+\left|\alpha\right|}{2}}\quad\forall x\in\R^n
\end{equation}
It follows from \eqref{Pr1Step2Eq4} and \eqref{Pr1Step2Eq5} that
\begin{align}\label{Pr1Step2Eq6}
\int_{B\(0,2r_0\)} z_\alpha\widetilde{U}_\mu\partial^{\(\alpha\)}\widetilde{U}_\mu\,dx&=\bigO\(\mu^{2k-\left|\alpha\right|}\int_{B\(0,2r_0/\mu\)}\big(1+\left|x\right|^2\big)^{-n+2k-\left|\alpha\right|/2}dx\)\nonumber\\
&=\left\{\begin{aligned}&\bigO\big(\mu^{2k-\left|\alpha\right|}\big)&&\text{if }\left|\alpha\right|>4k-n\\&\bigO\(\mu^{n-2k}\ln\(1/\mu\)\)&&\text{if }\left|\alpha\right|=4k-n\\&\bigO\(\mu^{n-2k}\)&&\text{if }\left|\alpha\right|<4k-n.\end{aligned}\right.
\end{align}
Finally, \eqref{Pr1Step2Eq1} follows from \eqref{Pr1Step2Eq2} and \eqref{Pr1Step2Eq6}. 
\endproof

We then prove the following:

\begin{step}\label{Pr1Step3}
Assume that $n\ge 2k+4$ and $g$ satisfies \eqref{Eq3} for some point $\xi\in M$. Then, as $\mu\to0$,
\begin{equation}\label{Pr1Step3Eq1}
\int_MU_\mu\Delta^kU_\mu\,dv_g=2^{2k-n}\(2k-1\)!\,\omega_n\Beta\(\frac{n}{2}-k,2k\)^{-1}+\bigO\(\mu^{n-2k}\).
\end{equation}
If $k\ge2$, then for every smooth function $f$ in $M$, 
\begin{align}\label{Pr1Step3Eq2}
&\int_MfU_\mu\Delta^{k-2}U_\mu \,dv_g=\frac{2^{2k-n-1}\(n-1\)!\(k-2\)!\,\omega_n}{\(n-2\)\(n-4\)\(n-2k-2\)}\,f\(\xi\)\mu^4\nonumber\\
&\qquad\qquad\times\sum_{l=k-2}^{2k-4}\frac{l!}{\(l-k+2\)!\(2k-l-4\)!\(n+l-2k-1\)!}\Beta\(\frac{n}{2}-k-1,l+1\)^{-1}\allowdisplaybreaks\nonumber\\
&\qquad\qquad\times\left\{\begin{aligned}&2\ln\(1/\mu\)&&\text{if }n=2k+4\text{ and }l=k-2\\&\Beta\(\frac{n}{2}+l-2k,2k-l-2\)&&\text{otherwise}\end{aligned}\right\}\nonumber\\
&\qquad\qquad+\left\{\begin{aligned}&\bigO\(\mu^4\)&&\text{if }n=2k+4\\&\smallo\(\mu^4\)&&\text{if }n>2k+4,\end{aligned}\right.
\end{align}
for every smooth, covariant tensor $T$ of rank $1$,   
\begin{align}\label{Pr1Step3Eq3}
&\int_M\(T,\nabla U_\mu\)\Delta^{k-2}U_\mu \,dv_g=-\frac{2^{2k-n-2}\(n-2k\)\(n-1\)!\(k-2\)!\,\omega_n}{\(n-2\)\(n-4\)\(n-2k-2\)}\tensor{T}{^a_{;a}}\(\xi\)\mu^4\nonumber\\
&\qquad\qquad\quad\times\sum_{l=k-2}^{2k-4}\frac{l!}{\(l-k+2\)!\(2k-l-4\)!\(n+l-2k\)!}\Beta\(\frac{n}{2}-k-1,l+1\)^{-1}\allowdisplaybreaks\nonumber\\
&\qquad\qquad\quad\times\left\{\begin{aligned}&2\ln\(1/\mu\)&&\text{if }n=2k+4\text{ and }l=k-2\\&\Beta\(\frac{n}{2}+l-2k,2k-l-2\)&&\text{otherwise}\end{aligned}\right\}\nonumber\\
&\qquad\qquad\quad+\left\{\begin{aligned}&\bigO\(\mu^4\)&&\text{if }n=2k+4\\&\smallo\(\mu^4\)&&\text{if }n>2k+4\end{aligned}\right.
\end{align}
and for every smooth, covariant tensor $T$ of rank $2$,   
\begin{align}\label{Pr1Step3Eq4}
&\int_M\(T,\nabla^2U_\mu\)\Delta^{k-2}U_\mu\,dv_g=\frac{2^{2k-n-4}\(n-2k\)\(n-1\)!\(k-2\)!\,\omega_n}{\(n-2\)\(n-4\)\(n-2k-2\)}\nonumber\\
&\qquad\qquad\times\sum_{l=k-2}^{2k-4}\frac{l!}{\(l-k+2\)!\(2k-l-4\)!\(n+l-2k+1\)!}\Beta\(\frac{n}{2}-k-1,l+1\)^{-1}\allowdisplaybreaks\nonumber\\
&\qquad\qquad\times\Bigg(-2\(n-4\)\(n+2l-2k\)\Beta\(\frac{n}{2}-2k+l+1,2k-l-2\)\tensor{T}{^a_{;a}}\(\xi\)\mu^2\allowdisplaybreaks\nonumber\\
&\qquad\qquad+\(\(n-2k+2\)\(\tensor{T}{^{ab}_{;ab}}\(\xi\)+\tensor{T}{^{ab}_{;ba}}\(\xi\)\)-\(n+2l-2k\)\tensor{T}{_a^a_{;b}^b}\(\xi\)\)\mu^4\allowdisplaybreaks\nonumber\\
&\qquad\qquad\times\left\{\begin{aligned}&2\ln\(1/\mu\)&&\text{if }n=2k+4\text{ and }l=k-2\\&\Beta\(\frac{n}{2}+l-2k,2k-l-2\)&&\text{otherwise}\end{aligned}\right\}\Bigg)\nonumber\\
&\qquad\qquad+\left\{\begin{aligned}&\bigO\(\mu^4\)&&\text{if }n=2k+4\\&\smallo\(\mu^4\)&&\text{if }n>2k+4.\end{aligned}\right.
\end{align}
If $k\ge3$, then for every smooth, covariant tensor $T$ of rank $2$,
\begin{align}\label{Pr1Step3Eq5}
&\int_M\(T,\nabla^2U_\mu\)\Delta^{k-3}U_\mu\,dv_g=-\frac{2^{2k-n-5}\(n-2k\)\(n-1\)!\(k-3\)!\,\omega_n}{\(n-2\)\(n-4\)\(n-2k-2\)}\tensor{T}{^a_{;a}}\(\xi\)\mu^4\nonumber\\
&\qquad\qquad\times\sum_{l=k-3}^{2k-6}\frac{\(n+2l-2k\)l!}{\(l-k+3\)!\(2k-l-6\)!\(n+l-2k+1\)!}\Beta\(\frac{n}{2}-k-1,l+1\)^{-1}\allowdisplaybreaks\nonumber\\
&\qquad\qquad\times\left\{\begin{aligned}&2\ln\(1/\mu\)&&\text{if }n=2k+4\text{ and }l=k-3\\&\Beta\(\frac{n}{2}+l-2k+1,2k-l-3\)&&\text{otherwise}\end{aligned}\right\}\nonumber\\
&\qquad\qquad+\left\{\begin{aligned}&\bigO\(\mu^4\)&&\text{if }n=2k+4\\&\smallo\(\mu^4\)&&\text{if }n>2k+4\end{aligned}\right.
\end{align}
and for every smooth, covariant tensor $T$ of rank $3$,
\begin{align}\label{Pr1Step3Eq6}
&\int_M\(T,\nabla^3U_\mu\)\Delta^{k-3}U_\mu\,dv_g=\frac{2^{2k-n-6}\(n-2k\)\(n-2k+2\)\(n-1\)!\(k-3\)!\,\omega_n}{\(n-2\)\(n-4\)\(n-2k-2\)}\nonumber\\
&\qquad\qquad\times\(\tensor{T}{_a^{ab}_{;b}}\(\xi\)+\tensor{T}{_a^{ba}_{;b}}\(\xi\)+\tensor{T}{^a_b^b_{;a}}\(\xi\)\)\mu^4\allowdisplaybreaks\nonumber\\
&\qquad\qquad\times\sum_{l=k-3}^{2k-6}\frac{\(n+2l-2k\)l!}{\(l-k+3\)!\(2k-l-6\)!\(n+l-2k+2\)!}\Beta\(\frac{n}{2}-k-1,l+1\)^{-1}\allowdisplaybreaks\nonumber\\
&\qquad\qquad\times\left\{\begin{aligned}&2\ln\(1/\mu\)&&\text{if }n=2k+4\text{ and }l=k-3\\&\Beta\(\frac{n}{2}+l-2k+1,2k-l-3\)&&\text{otherwise}\end{aligned}\right\}\nonumber\\
&\qquad\qquad+\left\{\begin{aligned}&\bigO\(\mu^4\)&&\text{if }n=2k+4\\&\smallo\(\mu^4\)&&\text{if }n>2k+4.\end{aligned}\right.
\end{align}
If $k\ge4$, then for every smooth, covariant tensor $T$ of rank $4$,
\begin{align}\label{Pr1Step3Eq7}
&\int_M\(T,\nabla^4U_\mu\)\Delta^{k-4}U_\mu\,dv_g=\frac{2^{2k-n-8}\(n-2k\)\(n-2k+2\)\(n-1\)!\(k-4\)!\,\omega_n}{3\(n-2\)\(n-4\)\(n-2k-2\)}\nonumber\\
&\qquad\qquad\times\(\tensor{T}{_a^a_b^b}\(\xi\)+\tensor{T}{_a^{ba}_b}\(\xi\)+\tensor{T}{_a^b_b^a}\(\xi\)\)\mu^4\nonumber\\
&\qquad\qquad\times\sum_{l=k-4}^{2k-8}\frac{\(n+2l-2k\)\(n+2l-2k+2\)l!}{\(l-k+4\)!\(2k-l-8\)!\(n+l-2k+3\)!}\Beta\(\frac{n}{2}-k-1,l+1\)^{-1}\allowdisplaybreaks\nonumber\\
&\qquad\qquad\times\left\{\begin{aligned}&2\ln\(1/\mu\)&&\text{if }n=2k+4\text{ and }l=k-4\\&\Beta\(\frac{n}{2}+l-2k+2,2k-l-4\)&&\text{otherwise}\end{aligned}\right\}\nonumber\\
&\qquad\qquad+\left\{\begin{aligned}&\bigO\(\mu^4\)&&\text{if }n=2k+4\\&\smallo\(\mu^4\)&&\text{if }n>2k+4.\end{aligned}\right.
\end{align}
\end{step}

\proof[Proof of Step~\ref{Pr1Step3}]
We let $j$ and $l$ be two integers such that 
$$\max\(2\(k-l-2\),0\)\le j\le k-l\quad\text{and}\quad\max\(k-4,0\)\le l\le k$$ 
and $T$ be a smooth, covariant tensor of rank $j$. By using geodesic normal coordinates, we obtain
\begin{align}\label{Pr1Step3Eq8}
&\int_M\(T,\nabla^jU_\mu\)\Delta^lU_\mu\,dv_g-\int_{B\(\xi,r_0\)}\(T,\nabla^jU_\mu\)\Delta^lU_\mu\,dv_g\nonumber\\
&\qquad=\int_{B\(0,2r_0\)\backslash B\(0,r_0\)}\widetilde{Z}_1\widetilde{U}_\mu\widetilde{Z}_2\widetilde{U}_\mu\,dx\nonumber\\
&\qquad=\sum_{\left|\alpha_1\right|\le j}\sum_{\left|\alpha_2\right|\le 2l}\int_{B\(0,2r_0\)\backslash B\(0,r_0\)} z_{1,\alpha}z_{2,\alpha}\big(\partial^{\alpha_2}\widetilde{U}_\mu\big)\big(\partial^{\alpha_1}\widetilde{U}_\mu\big)dx,
\end{align}
where $\widetilde{U}_\mu$ is as in \eqref{Pr1Step2Eq3} and
$$\widetilde{Z}_1\(x\):=\sum_{\left|\alpha\right|\le j}z_{1,\alpha}\(x\)\partial^{\alpha}\quad\text{and}\quad\widetilde{Z}_2\(x\):=\sum_{\left|\alpha\right|\le2l}z_{2,\alpha}\(x\)\partial^{\alpha}\quad\forall x\in B\(0,2r_0\)$$
for some smooth functions $z_{1,\alpha}$ and $z_{2,\alpha}$ in $B\(0,2r_0\)$. By proceeding as in \eqref{Pr1Step2Eq4}--\eqref{Pr1Step2Eq6}, we obtain
\begin{equation}\label{Pr1Step3Eq9}
\int_{B\(0,2r_0\)\backslash B\(0,r_0\)} z_{1,\alpha}z_{2,\alpha}\partial^{\alpha_2}\widetilde{U}_\mu\partial^{\alpha_1}\widetilde{U}_\mu\,dx=\bigO\(\mu^{n-2k}\).
\end{equation}
It follows from \eqref{Pr1Step3Eq8} and \eqref{Pr1Step3Eq9} that
\begin{equation}\label{Pr1Step3Eq10}
\int_M\(T,\nabla^jU_\mu\)\Delta^lU_\mu\,dv_g=\int_{B\(\xi,r_0\)}\(T,\nabla^jU_\mu\)\Delta^lU_\mu\,dv_g+\bigO\(\mu^{n-2k}\).
\end{equation}
By using \eqref{Eq3} and rewriting the integral in the right-hand side of \eqref{Pr1Step3Eq10} in geodesic normal coordinates, we obtain
\begin{equation}\label{Pr1Step3Eq11}
\int_{B\(\xi,r_0\)}\(T,\nabla^jU_\mu\)\Delta^lU_\mu\,dv_g=\sum_{j'=0}^j\int_{B\(0,r_0\)}\hspace{-1pt}\widehat{T}^{i_1\dotsc i_{j'}}\circ\exp_\xi U_{\mu,i_1\dotsc i_{j'}}\Delta_0^l\,U_\mu\,dx,
\end{equation}
where $U_{\mu,i_1\dotsc i_{j'}}:=\partial^{(i_1\dotsc i_{j'})}\(U_\mu\circ\,\exp_\xi\)$ and the tensor $\widehat{T}$ is defined as
$$\widehat{T}^{i_1\dotsc i_{j'}}:=\left\{\begin{aligned}&T^{i_1\dotsc i_j}&&\text{if }j'=j\\&-\Gamma^{i_1\dotsc i_{j'}}_{e_1\dotsc e_j}T^{e_1\dotsc e_j}&&\text{if }j'<j\end{aligned}\right\},$$
where $\Gamma^{i_1\dotsc i_{j'}}_{e_1\dotsc e_j}$ is the generalized Christoffel symbol such that $\Gamma^{i_1\dotsc i_{j'}}_{e_1\dotsc e_j}$ is symmetric in $i_1,\dotsc,i_{j'}$ and
$$u_{;e_1\dotsc e_j}=u_{,e_1\dotsc e_j}-\sum_{j'=0}^{j-1}\Gamma^{i_1\dotsc i_{j'}}_{e_1\dotsc e_j}u_{,i_1\dotsc i_{j'}}$$
in geodesic normal coordinates. By using \eqref{Pr1Step3Eq11} together with a  straightforward change of variable and a Taylor expansion, we then obtain
\begin{align}\label{Pr1Step3Eq12}
&\int_{B\(\xi,r_0\)}\(T,\nabla^jU_\mu\)\Delta^lU_\mu\,dv_g\nonumber\\
&\qquad\qquad=\sum_{j'=0}^j\mu^{2k-2l-j'}\int_{B\(0,r_0/\mu\)}\widehat{T}^{i_1\dotsc i_{j'}}\(\exp_\xi\(\mu x\)\)U_{,i_1\dotsc i_{j'}}\(x\)\Delta_0^l\,U\(x\)dx\allowdisplaybreaks\nonumber\\
&\qquad\qquad=\sum_{j'=\max\(2\(k-l-2\),0\)}^j\sum_{j''=0}^{j'+2l-2k+4}\frac{\mu^{2k-2l-j'+j''}}{j''!}\widehat{T}^{i_1\dotsc i_{j'},i_{j'+1}\dotsc i_{j'+j''}}\(\xi\)\nonumber\\
&\qquad\qquad\quad\times\int_{B\(0,r_0/\mu\)}U_{,i_1\dotsc i_{j'}}x_{i_{j'+1}}\dotsm x_{i_{j'+j''}}\Delta_0^l\,U\,dx+\bigO\Bigg(\sum_{j'=0}^j\mu^{\max(5,2k-2l-j')}\nonumber\\
&\qquad\qquad\quad\times\int_{B\(0,r_0/\mu\)}\left|x\right|^{\max(j'+2l-2k+5,0)}\big|U_{,i_1\dotsc i_{j'}}\Delta_0^l\,U\big|\,dx\Bigg).
\end{align}
On the other hand, by using \eqref{Pr1Step2Eq5}, we obtain
\begin{align}\label{Pr1Step3Eq13}
&\mu^{\max(5,2k-2l-j')}\int_{B\(0,r_0/\mu\)}\left|x\right|^{\max(j'+2l-2k+5,0)}\big|U_{,i_1\dotsc i_{j'}}\Delta_0^l\,U\big|\,dx\nonumber\\
&\quad=\bigO\(\mu^{\max(5,2k-2l-j')}\int_{B\(0,r_0/\mu\)}\left|x\right|^{\max(j'+2l-2k+5,0)}\big(1+\left|x\right|^2\big)^{-\frac{2n+j'+2l-4k}{2}}dx\)\nonumber\\
&\quad=\left\{\begin{aligned}&\bigO\(\mu^4\)&&\text{if }n=2k+4\\&\smallo\(\mu^4\)&&\text{otherwise.}\end{aligned}\right.
\end{align}
It follows from \eqref{Pr1Step3Eq10}, \eqref{Pr1Step3Eq12} and \eqref{Pr1Step3Eq13} that
\begin{align}\label{Pr1Step3Eq14}
&\int_M\(T,\nabla^jU_\mu\)\Delta^lU_\mu\,dv_g=\sum_{j'=\max\(2\(k-l-2\),0\)}^{j}\sum_{j''=0}^{j'+2l-2k+4}\frac{\mu^{2k-2l-j'+j''}}{j''!}\nonumber\\
&\qquad\qquad\qquad\quad\times\widehat{T}^{i_1\dotsc i_{j'},i_{j'+1}\dotsc i_{j'+j''}}\(\xi\)\int_{B\(0,r_0/\mu\)}U_{,i_1\dotsc i_{j'}}x_{i_{j'+1}}\dotsm x_{i_{j'+j''}}\Delta_0^l\,U\,dx\nonumber\\
&\qquad\qquad\qquad\quad+\left\{\begin{aligned}&\bigO\(\mu^4\)&&\text{if }n=2k+4\\&\smallo\(\mu^4\)&&\text{if }n>2k+4.\end{aligned}\right.
\end{align}
An easy induction gives
\begin{multline}\label{Pr1Step3Eq15}
U_{,i_1\dotsc i_j}\(x\)=\sum_{m=0}^{\lfloor j/2\rfloor}\frac{2^{j-2m}}{m!\(j-2m\)!}\,\partial_r^{j-m}U\(r\)\\
\times\sum_{\sigma\in\mathfrak{S}(j)}\delta_{i_{\sigma\(1\)}i_{\sigma\(2\)}}\dotsm\delta_{i_{\sigma\(2m-1\)}i_{\sigma\(2m\)}}x_{i_{\sigma\(2m+1\)}}\dotsm x_{i_{\sigma\(j\)}}\quad\forall x\in\R^n,
\end{multline}
where $r:=\left|x\right|^2$, $U\(r\):=U\(x\)=\(1+r\)^{\(2k-n\)/2}$, $\mathfrak{S}\(j\)$ is the set of all permutations of $\(1,\dotsc,j\)$ and $\delta_{i_{\sigma\(1\)}i_{\sigma\(2\)}},\dotsc,\delta_{i_{\sigma\(2m-1\)}i_{\sigma\(2m\)}}$ stand for the Kronecker symbols. Furthermore, it is easy to see that
\begin{align}\label{Pr1Step3Eq16}
\partial_r^{j}U\(r\)&=\(-1\)^j2^{-j}\(n-2k\)\(n-2k+2\)\dotsm\(n-2k+2j-2\)\(1+r\)^{-\frac{n-2k+2j}{2}}\nonumber\allowdisplaybreaks\\
&=\frac{2\(-1\)^jj!}{\(n-2k-2\)}\Beta\(\frac{n}{2}-k-1,j+1\)^{-1}\(1+r\)^{-\frac{n-2k+2j}{2}}.
\end{align}
Another induction yields
\begin{equation}\label{Pr1Step3Eq17}
\Delta_0^l\,U\(x\)=\left\{\begin{aligned}&\frac{2^{2l+1}\,l!}{\(n-2k-2\)\(k-l-1\)!}\sum_{l'=l}^{2l}\frac{l'!\(k+l-l'-1\)!}{\(l'-l\)!\(2l-l'\)!}\\
&\qquad\times\Beta\(\frac{n}{2}-k-1,l'+1\)^{-1}\(1+r\)^{-\frac{n+2l'-2k}{2}}&&\text{if }l<k\\&2^{2k}\(2k-1\)!\Beta\(\frac{n}{2}-k,2k\)^{-1}\(1+r\)^{-\frac{n+2k}{2}}&&\text{if }l=k\end{aligned}\right.
\end{equation}
for all $x\in\R^n$. In the case where $j=0$, $l=k$ and $T\equiv1$, it follows from \eqref{Pr1Step3Eq17} that
\begin{equation}\label{Pr1Step3Eq18}
\int_{B\(0,r_0/\mu\)}U\Delta^k_0\,U\,dx=2^{2k-1}\(2k-1\)!\,\omega_{n-1}\Beta\(\frac{n}{2}-k,2k\)^{-1}\int_0^{\(r_0/\mu\)^2}\frac{r^{\frac{n-2}{2}}}{\(1+r\)^n}\,dr,
\end{equation}
where $\omega_{n-1}=\Vol\(\S^{n-1},g_0\)$ is the volume of the standard $\(n-1\)$-dimensional sphere. On the other hand, in the case where $l<k$, by putting together \eqref{Pr1Step3Eq15}--\eqref{Pr1Step3Eq17} , we obtain 
\begin{align}\label{Pr1Step3Eq19}
&\int_{B\(0,r_0/\mu\)}U_{,i_1\dotsc i_{j'}}x_{i_{j'+1}}\dotsm x_{i_{j'+j''}}\Delta^l_0\,U\,dx=\frac{2^{2l+1}\,l!}{\(n-2k-2\)^2\(k-l-1\)!}\sum_{l'=l}^{2l}\sum_{m=0}^{\lfloor j'/2\rfloor}\nonumber\\
&\ \ \ \frac{\(-1\)^{j'-m}2^{j'-2m}\,l'!\(k+l-l'-1\)!\(j'-m\)!}{\(l'-l\)!\(2l-l'\)!\,m!\(j'-2m\)!}\Beta\(\frac{n}{2}-k-1,l'+1\)^{-1}\allowdisplaybreaks\nonumber\\
&\ \ \ \times\Beta\(\frac{n}{2}-k-1,j'-m+1\)^{-1}\int_0^{\(r_0/\mu\)^2}\frac{r^{\frac{n+j'+j''-2m-2}{2}}}{\(1+r\)^{n+j'-m+l'-2k}}\,dr\sum_{\sigma\in\mathfrak{S}(j')}\delta_{i_{\sigma(1)}i_{\sigma(2)}}\nonumber\\
&\ \ \ \dotsm\delta_{i_{\sigma(2m-1)}i_{\sigma(2m)}}\int_{\S^{n-1}}y_{i_{\sigma(2m+1)}}\dotsm y_{i_{\sigma(j')}}y_{i_{j'+1}}\dotsm y_{i_{j'+j''}}dv_{g_0}\(y\).
\end{align}
A standard computation gives
\begin{equation}\label{Pr1Step3Eq20}
\int_0^{\(r_0/\mu\)^2}\frac{r^{a-1}dr}{\(1+r\)^b}=\left\{\begin{aligned}&2\ln\(1/\mu\)+\bigO\(1\)&&\text{if }b=a\\&\Beta\(a,b-a\)+\bigO\big(\mu^{2\(b-a\)}\big)&&\text{if }b>a.\end{aligned}\right.
\end{equation}
On the other hand, by using the fact (see for example Brendle~\cite{Bre}*{Proposition~28}) that for every homogeneous polynomial $\varPhi$ of degree $j\ge2$, 
$$\int_{\S^{n-1}}\varPhi\(y\)\,dv_{g_0}\(y\)=\frac{-1}{j\(n+j-2\)}\int_{\S^{n-1}}\Delta_{0}\varPhi\(y\)\,dv_{g_0}\(y\),$$ 
another induction yields that when $j$ is even,
\begin{multline}\label{Pr1Step3Eq21}
\int_{\S^{n-1}}y_{i_1}\dotsm y_{i_j}\,dv_{g_0}\(y\)=\frac{\(n-2\)\omega_{n-1}}{2^{j+1}\(j/2\)!^2}\Beta\(\frac{n-2}{2},\frac{j+2}{2}\)\\
\times\sum_{\sigma\in\mathfrak{S}(j)}\delta_{i_{\sigma\(1\)}i_{\sigma\(2\)}}\dotsm\delta_{i_{\sigma\(j-1\)}i_{\sigma\(j\)}}.
\end{multline}
The integral in \eqref{Pr1Step3Eq21} vanishes when $j$ is odd. By observing that 
\begin{equation}\label{Pr1Step3Eq22}
\omega_n=2^{n-1}\Beta\(\frac{n}{2},\frac{n}{2}\)\omega_{n-1},
\end{equation}
we obtain that for even $j$,
\begin{equation}\label{Pr1Step3Eq23}
\Beta\(\frac{n-2}{2},\frac{j+2}{2}\)=\frac{2^{2-n}\(n-1\)!\(j/2\)!\,\omega_n}{\(n-2\)\(n+j/2-1\)!\,\omega_{n-1}}\Beta\(\frac{n}{2},\frac{n+j}{2}\)^{-1}.
\end{equation}
By using \eqref{Pr1Step3Eq20}--\eqref{Pr1Step3Eq23} together with the identity
\begin{align*}
&\Beta\(\frac{n}{2},\frac{n+j'+j''-2m}{2}\)^{-1}\Beta\(\frac{n+j'+j''-2m}{2},\frac{n+j'-j''+2l'-4k}{2}\)\\
&\qquad\qquad\qquad\qquad\qquad=\frac{\(\frac{j'+j''}{2}+n-m-1\)!}{\(n+j'-m+l'-2k-1\)!\(\frac{j''-j'}{2}+2k-l'-1\)!}\\
&\qquad\qquad\qquad\qquad\qquad\quad\times\Beta\(\frac{n+j'-j''+2l'-4k}{2},\frac{j''-j'+4k-2l'}{2}\),
\end{align*}
we obtain that if $j'+j''$ is even, then
\begin{align}\label{Pr1Step3Eq24}
&\int_0^{\(r_0/\mu\)^2}\frac{r^{\frac{n+j'+j''-2m-2}{2}}}{\(1+r\)^{n+j'-m+l'-2k}}\,dr\int_{\S^{n-1}}y_{i_{\sigma(2m+1)}}\dotsm y_{i_{\sigma(j')}}y_{i_{j'+1}}\dotsm y_{i_{j'+j''}}dv_{g_0}\(y\)\nonumber\\
&=\frac{2^{1-n-j'-j''+2m}\(n-1\)!\,\omega_n}{\(n+j'-m+l'-2k-1\)!\(\frac{j''-j'}{2}+2k-l'-1\)!\(\frac{j'+j''}{2}-m\)!}\nonumber\\
&\times\left\{\begin{aligned}&2\ln\(1/\mu\)+\bigO\(1\)\hspace{117pt}\text{if }n+j'-j''+2l'-4k=0\\&\Beta\(\frac{n+j'-j''+2l'-4k}{2},\frac{j''-j'+4k-2l'}{2}\)\nonumber\\
&\hspace{80pt}+\bigO\big(\mu^{n+j'-j''+2l'-4k}\big)\quad\text{if }0<n+j'-j''+2l'-4k<n\end{aligned}\right\}\nonumber\\
&\times\sum_{\sigma'\in\mathfrak{S}(S_{j',j'',m,\sigma})}\delta_{i_{\sigma'(\sigma(2m+1))}i_{\sigma'(\sigma(2m+2))}}\dotsm\delta_{i_{\sigma'(\sigma(j'-1))}i_{\sigma'(\sigma(j'))}}\nonumber\\
&\times\delta_{i_{\sigma'(j'+1)}i_{\sigma'(j'+2)}}\dotsm\delta_{i_{\sigma'(j'+j''-1)}i_{\sigma'(j'+j'')}},
\end{align} 
where 
$$S_{j',j'',m,\sigma}:=\big(\sigma\(2m+1\),\dotsc,\sigma\(j'\),j'+1,\dotsc,j'+j''\big)$$
and $\mathfrak{S}\(S_{j',j'',m,\sigma}\)$ stands for the set of all permutations of $S_{j',j'',m,\sigma}$. In the case where $j=0$, $l=k$ and $T\equiv1$, \eqref{Pr1Step3Eq1} follows from \eqref{Pr1Step3Eq10}, \eqref{Pr1Step3Eq11}, \eqref{Pr1Step3Eq18}, \eqref{Pr1Step3Eq20} and \eqref{Pr1Step3Eq22}. On the other hand, in the case where $l<k$, by combining \eqref{Pr1Step3Eq14}, \eqref{Pr1Step3Eq19} and \eqref{Pr1Step3Eq24} (and replacing $j''$ by $j'-2m'+2l-2k+4$ for $m'\in\left\{0,\dotsc,\lfloor j'/2\rfloor+l-k+2\right\}$ so that $j'+j''$ is even and $0\le j''\le j'+2l-2k+4$), we obtain
\begin{align}\label{Pr1Step3Eq25}
&\int_M\(T,\nabla^jU_\mu\)\Delta^lU_\mu\,dv_g=\frac{2^{2k-n-2}\(n-1\)!\,l!\,\omega_n}{\(n-2k-2\)^2\(k-l-1\)!}\sum_{l'=l}^{2l}\sum_{j'=\max\(2\(k-l-2\),0\)}^{j}\sum_{m=0}^{\lfloor j'/2\rfloor}\nonumber\\
&\sum_{m'=0}^{\lfloor j'/2\rfloor+l-k+2}\frac{2^{2m'-j'}l'!\(k+l-l'-1\)!\,c\(n,k,j',l,l',m,m'\)\mu^{4-2m'}}{\(l'-l\)!\(2l-l'\)!\(k+l-l'-m'+1\)!\(j'-2m'+2l-2k+4\)!}\allowdisplaybreaks\nonumber\\
&\times\Beta\(\frac{n}{2}-k-1,l'+1\)^{-1}\widehat{T}^{i_1\dotsc i_{j'},i_{j'+1}\dotsc i_{2(j'-m'+l-k+2)}}\(\xi\)\allowdisplaybreaks\nonumber\\
&\times\sum_{\sigma\in\mathfrak{S}(j')}\sum_{\sigma'\in\mathfrak{S}(S_{j',j'-2m'+2l-2k+4,m,\sigma})}\delta_{i_{\sigma(1)}i_{\sigma(2)}}\dotsm\delta_{i_{\sigma(2m-1)}i_{\sigma(2m)}}\allowdisplaybreaks\nonumber\\
&\times\delta_{i_{\sigma'(\sigma(2m+1))}i_{\sigma'(\sigma(2m+2))}}\dotsm\delta_{i_{\sigma'(\sigma(j'-1))}i_{\sigma'(\sigma(j'))}}\allowdisplaybreaks\nonumber\\
&\times\delta_{i_{\sigma'(j'+1)}i_{\sigma'(j'+2)}}\dotsm\delta_{i_{\sigma'(2(j'-m'+l-k+2)-1)}i_{\sigma'(2(j'-m'+l-k+2))}}\allowdisplaybreaks\nonumber\\
&\times\left\{\begin{aligned}&2\ln\(1/\mu\)+\bigO\(1\)\hspace{115pt}\text{if }n=2k+4,\ l'=l\text{ and }m'=0\\&\Beta\(\frac{n}{2}+m'+l'-l-k-2,k+l-l'-m'+2\)+\smallo\big(\mu^{2m'}\big)\hspace{33pt}\text{otherwise}\end{aligned}\right\}\nonumber\\
&+\left\{\begin{aligned}&\bigO\(\mu^4\)&&\text{if }n=2k+4\\&\smallo\(\mu^4\)&&\text{if }n>2k+4,\end{aligned}\right.
\end{align}
where
\begin{multline*}
c\(n,k,j',l,l',m,m'\):=\(-1\)^{j'-m}\Beta\(\frac{n}{2}-k-1,j'-m+1\)^{-1}\\
\times\frac{\(j'-m\)!}{m!\(j'-2m\)!\(n+j'-m+l'-2k-1\)!\(j'-m-m'+l-k+2\)!}\,.
\end{multline*}
Straightforward computations yield
\begin{align}\label{Pr1Step3Eq26}
&\sum_{\sigma\in\mathfrak{S}(j')}\sum_{\sigma'\in\mathfrak{S}(S_{j',j'',m,\sigma})}\delta_{i_{\sigma(1)}i_{\sigma(2)}}\dotsm\delta_{i_{\sigma(2m-1)}i_{\sigma(2m)}}\delta_{i_{\sigma'(\sigma(2m+1))}i_{\sigma'(\sigma(2m+2))}}\allowdisplaybreaks\nonumber\\
&\qquad\dotsm\delta_{i_{\sigma'(\sigma(j'-1))}i_{\sigma'(\sigma(j'))}}\delta_{i_{\sigma'(j'+1)}i_{\sigma'(j'+2)}}\dotsm\delta_{i_{\sigma'(j'+j''-1)}i_{\sigma'(j'+j'')}}\allowdisplaybreaks\nonumber\\
&\quad=\left\{\begin{aligned}
&1&&\text{if }j'=j''=m=0\\
&2\,\delta_{i_1i_2}&&\text{if }j'=j''=1\text{ and }m=0\\
&2\(2-m\)\delta_{i_1i_2}&&\text{if }j'=2,\ j''=0\text{ and }m\le1\\
&16\(\delta_{i_1i_2}\delta_{i_3i_4}+\delta_{i_1i_3}\delta_{i_2i_4}+\delta_{i_1i_4}\delta_{i_2i_3}\)&&\text{if }j'=j''=2\text{ and }m=0\\
&4\,\delta_{i_1i_2}\delta_{i_3i_4}&&\text{if }j'=j''=2\text{ and }m=1\\
&2\(4-2m\)!\(\delta_{i_1i_2}\delta_{i_3i_4}+\delta_{i_1i_3}\delta_{i_2i_4}+\delta_{i_1i_4}\delta_{i_2i_3}\)&&\text{if }j'=3,\ j''=1\text{ and }m\le1\\
&8\(4-2m\)!\(\delta_{i_1i_2}\delta_{i_3i_4}+\delta_{i_1i_3}\delta_{i_2i_4}+\delta_{i_1i_4}\delta_{i_2i_3}\)&&\text{if }j'=4,\ j''=0\text{ and }m\le2.\\
\end{aligned}\right.
\end{align}
On the other hand, by using \eqref{Eq4} and the fact that for all $a,b,c,d,e\in\left\{1,\dotsc,n\right\}$, in geodesic normal coordinates,
$$g_{ab}\(\xi\)=\delta_{ab},\quad g_{ab,c}\(\xi\)=0\quad\text{and}\quad g_{ab,cd}\(\xi\)=\frac{1}{3}\(\Riemann_{acdb}\(\xi\)+\Riemann_{adcb}\(\xi\)\),$$
and
$$\Gamma^a_{bc}\(\xi\)=0,\quad\Gamma^a_{bc,d}\(\xi\)=\frac{1}{3}\(\Riemann_{abdc}\(\xi\)+\Riemann_{acdb}\(\xi\)\)\quad\text{and}\quad\Gamma^{ab}_{cde}\(\xi\)=0,$$
we obtain
\begin{equation}\label{Pr1Step3Eq27}
\left\{\begin{aligned}
&\widehat{T}\(\xi\)=T\(\xi\)&&\text{if }j=0\\
&\tensor{\widehat{T}}{^a_{,a}}\(\xi\)=\tensor{T}{^a_{;a}}\(\xi\)&&\text{if }j=1\\
&\left\{\begin{aligned}
&\tensor{\widehat{T}}{^a_{,a}}\(\xi\)=0,\quad\tensor{\widehat{T}}{_a^a}\(\xi\)=\tensor{T}{_a^a}\(\xi\),\quad\tensor{\widehat{T}}{_a^a_{,b}^b}\(\xi\)=\tensor{T}{_a^a_{;b}^b}\(\xi\)\\
&\hspace{77pt}\text{and}\quad\tensor{\widehat{T}}{^{ab}_{,ab}}\(\xi\)=\tensor{\widehat{T}}{^{ab}_{,ba}}\(\xi\)=\tensor{T}{^{ab}_{;ab}}\(\xi\)
\end{aligned}\right\}&&\text{if }j=2\\
&\left\{\begin{aligned}
&\tensor{\widehat{T}}{_a^a}\(\xi\)=0\quad\text{and}\quad\tensor{\widehat{T}}{_a^{ab}_{,b}}\(\xi\)+\tensor{\widehat{T}}{_a^{ba}_{,b}}\(\xi\)+\tensor{\widehat{T}}{^a_b^b_{,a}}\(\xi\)\\
&\hspace{93pt}=\tensor{T}{_a^{ab}_{;b}}\(\xi\)+\tensor{T}{_a^{ba}_{;b}}\(\xi\)+\tensor{T}{^a_b^b_{;a}}\(\xi\)
\end{aligned}\right\}&&\text{if }j=3\\
&\tensor{\widehat{T}}{_a^a_b^b}\(\xi\)+\tensor{\widehat{T}}{_{ab}^{ab}}\(\xi\)+\tensor{\widehat{T}}{_{ab}^{ba}}\(\xi\)=\tensor{T}{_a^a_b^b}\(\xi\)+\tensor{T}{_{ab}^{ab}}\(\xi\)+\tensor{T}{_{ab}^{ba}}\(\xi\)&&\text{if }j=4.
\end{aligned}\right.
\end{equation}
We then obtain \eqref{Pr1Step3Eq2} by putting together \eqref{Pr1Step3Eq25}, \eqref{Pr1Step3Eq26} and \eqref{Pr1Step3Eq27} and using the identities
$$c\(n,k,0,k-2,l',0,0\)=\frac{n-2k-2}{2\(n+l'-2k-1\)!}$$
and
\begin{multline*}
\Beta\(\frac{n}{2}+l'-2k,2k-l'+4\)\\
=\frac{4\(2k-l'-1\)\(2k-l'-2\)}{\(n-2\)\(n-4\)}\Beta\(\frac{n}{2}+l'-2k,2k-l'-2\).
\end{multline*}
The estimates \eqref{Pr1Step3Eq3}--\eqref{Pr1Step3Eq7} follow in the same way from \eqref{Pr1Step3Eq25}, \eqref{Pr1Step3Eq26} and \eqref{Pr1Step3Eq27} by using the identities
\begin{align*}
&c\(n,k,1,k-2,l',0,0\)=-\frac{\(n-2k-2\)\(n-2k\)}{4\(n+l'-2k\)!}\,,\allowdisplaybreaks\\
&c\(n,k,2,k-2,l',0,0\)=\frac{\(n-2k-2\)\(n-2k\)\(n-2k+2\)}{32\(n+l'-2k+1\)!}\,,\allowdisplaybreaks\\
&2c\(n,k,2,k-2,l',0,1\)+c\(n,k,2,k-2,l',1,1\)\\
&\qquad=4c\(n,k,2,k-2,l',0,0\)+c\(n,k,2,k-2,l',1,0\)\allowdisplaybreaks\\
&\qquad=2c\(n,k,2,k-3,l',0,0\)+c\(n,k,2,k-3,l',1,0\)\\
&\qquad=-\frac{\(n-2k-2\)\(n-2k\)\(n+2l'-2k\)}{8\(n+l'-2k+1\)!}\,,\allowdisplaybreaks\\
&24c\(n,k,3,k-3,l',0,0\)+2c\(n,k,3,k-3,l',1,0\)\allowdisplaybreaks\\
&\qquad=\frac{\(n-2k-2\)\(n-2k\)\(n-2k+2\)\(n+2l'-2k\)}{8\(n+l'-2k+2\)!},\allowdisplaybreaks\\
&24c\(n,k,4,k-4,l',0,0\)+2c\(n,k,4,k-4,l',1,0\)+c\(n,k,4,k-4,l',2,0\)\\
&\qquad=\frac{\(n-2k-2\)\(n-2k\)\(n-2k+2\)\(n+2l'-2k\)\(n+2l'-2k+2\)}{64\(n+l'-2k+3\)!}
\end{align*}
and
\begin{align*}
&\Beta\(\frac{n}{2}+l'-l-k-2,k+l-l'+2\)\nonumber\\
&\qquad=\frac{2\(k+l-l'+1\)}{n-2}\Beta\(\frac{n}{2}+l'-l-k-2,k+l-l'+1\)\nonumber\allowdisplaybreaks\\
&\qquad=\frac{4\(k+l-l'+1\)\(k+l-l'\)}{\(n-2\)\(n-4\)}\Beta\(\frac{n}{2}+l'-l-k-2,k+l-l'\).
\end{align*}
This ends the proof of Step~\ref{Pr1Step3}.
\endproof

As regards the integral in the denominator of $I_{k,f,g}\(u\)$, we obtain the following:

\begin{step}\label{Pr1Step4}
Assume that $n\ge 2k+1$ and $g$ satisfies \eqref{Eq3} for some point $\xi\in M$. Then, for every smooth function $f$ in $M$, as $\mu\to0$,
\begin{equation}\label{Pr1Step4Eq}
\int_MfU_\mu^{2^*_k}\,dv_g=\frac{\omega_n}{2^n}f\(\xi\)-\frac{\omega_n\Delta f\(\xi\)\mu^2}{2^{n+1}\(n-2\)}\\+\frac{\omega_n\Delta^2f\(\xi\)\mu^4}{2^{n+3}\(n-2\)\(n-4\)}+\smallo\(\mu^4\).
\end{equation}
\end{step}

\proof[Proof of Step~\ref{Pr1Step4}]
By observing that $U_\mu^{2^*_k}$ does not depend on $k$ in $B\(0,r_0\)$, we obtain that \eqref{Pr1Step4Eq} is in fact identical to an estimate obtained by Esposito and Robert~\cite{EspRob} in the case where $k=2$ (note that in our case, $\Ricci\(\xi\)=0$ and $\nabla\Scal\(\xi\)=0$ since we are working with conformal normal coordinates, see \eqref{Eq4} and \eqref{Eq5}). 
\endproof

We can now end the proof of Proposition~\ref{Pr1} by putting together the results of Steps~\ref{Pr1Step1}--~\ref{Pr1Step4}:

\proof[End of proof of Proposition~\ref{Pr1}]
We assume that $k\ge2$ and refer to Aubin~\cite{Aub} for the case where $k=1$. By using \eqref{Pr1Step4Eq}, we obtain
\begin{multline}\label{Pr1Eq2}
\(\int_MfU_\mu^{2^*_k}\,dv_g\)^{-\frac{n-2k}{n}}=\(\frac{\omega_n}{2^n}f\(\xi\)\)^{-\frac{n-2k}{n}}
\Bigg[1+\frac{\(n-2k\)\Delta f\(\xi\)\mu^2}{2n\(n-2\)f\(\xi\)}\\
-\frac{n-2k}{4n\(n-2\)}\(\frac{\Delta^2f\(\xi\)}{2\(n-4\)f\(\xi\)}-\frac{\(n-k\)\(\Delta f\(\xi\)\)^2}{n\(n-2\)f\(\xi\)^2}\)\mu^4+\smallo\(\mu^4\)\Bigg].
\end{multline}
We let $J_1$, $J_2$,  $T_1$, $T_2$, $T_3$, $T_4$, $T_5$ and $Z$ be as in Step~\ref{Pr1Step1}. 
Since $k\ge1$, by integrating by parts, we obtain
\begin{multline}\label{Pr1Eq3}
\int_MU_\mu\Delta^{k-1}\(J_1U_\mu\)dv_g=\int_M\Delta\(J_1U_\mu\)\Delta^{k-2}U_\mu\,dv_g\\
=\int_M\(U_\mu\Delta J_1-2\(\nabla J_1,\nabla U_\mu\)+J_1\Delta U_\mu\)\Delta^{k-2}U_\mu\,dv_g.
\end{multline}
By integrating by parts again, it follows from \eqref{Pr1Step1Eq1} and \eqref{Pr1Eq3} that
\begin{align}\label{Pr1Eq4}
&\int_MU_\mu P_{2k}U_\mu\,dv_g=\int_MU_\mu\Delta^kU_\mu\,dv_g+k\int_M\big(\(\(k-1\)J_2+\Delta J_1\)U_\mu\nonumber\\
&\quad+\(\(k-1\)T_1-2\nabla J_1,\nabla U_\mu\)+\(\(k-1\)T_2-J_1g,\nabla^2U_\mu\)\big)\Delta^{k-2}U_\mu\,dv_g\allowdisplaybreaks\nonumber\\
&\quad+k\(k-1\)\(k-2\)\int_M\(\(T_3,\nabla^2U_\mu\)+\(T_4,\nabla^3U_\mu\)\)\Delta^{k-3}U_\mu\,dv_g\nonumber\\
&\quad+k\(k-1\)\(k-2\)\(k-3\)\int_M\(T_5,\nabla^4U_\mu\)\Delta^{k-4}U_\mu dv_g+\int_MU_\mu ZU_\mu\,dv_g.
\end{align}
By using \eqref{Eq5}, we obtain 
\begin{equation}\label{Pr1Eq5}
\tensor{\Schouten}{_a^a_{;b}^b}\(\xi\)=\tensor{\Schouten}{^a_{b;a}^b}\(\xi\)=\tensor{\Schouten}{_a^b_{;b}^a}\(\xi\)=-\frac{\left|\Weyl\(\xi\)\right|^2}{12\(n-1\)}\,.
\end{equation}
By using \eqref{Eq4}, \eqref{Eq5} and \eqref{Pr1Eq5} together with straightforward computations, we obtain
$$\tensor{(T_3)}{_a^a}\(\xi\)=-\frac{n+3k+1}{36\(n-1\)}\left|\Weyl\(\xi\)\right|^2=-\frac{n+3k+1}{36n\(n-1\)}\left|\Weyl\(\xi\)\right|^2\tensor{g}{_a^a}$$
and
\begin{align*}
&\tensor{(T_4)}{_a^{ab}_{;b}}\(\xi\)+\tensor{(T_4)}{_a^{ba}_{;b}}\(\xi\)+\tensor{(T_4)}{^a_b^b_{;a}}\(\xi\)=-\frac{k+1}{6\(n-1\)}\left|\Weyl\(\xi\)\right|^2\\
&\qquad=\frac{k+1}{\(n-1\)\(n+2\)}\big(\tensor{(\nabla S\otimes g)}{_a^{ab}_{;b}}\(\xi\)+\tensor{(\nabla S\otimes g)}{_a^{ba}_{;b}}\(\xi\)+\tensor{(\nabla S\otimes g)}{^a_b^b_{;a}}\(\xi\)\big)
\end{align*}
and
\begin{multline*}
\tensor{(T_5)}{_a^a_b^b}\(\xi\)+\tensor{(T_5)}{_{ab}^{ab}}\(\xi\)+\tensor{(T_5)}{_{ab}^{ba}}\(\xi\)=-\frac{k+1}{10\(n-1\)}\left|\Weyl\(\xi\)\right|^2\\
=-\frac{\(k+1\)\left|\Weyl\(\xi\)\right|^2}{10n\(n-1\)\(n+2\)}\big(\tensor{(g\otimes g)}{_a^a_b^b}\(\xi\)+\tensor{(g\otimes g)}{_{ab}^{ab}}\(\xi\)+\tensor{(g\otimes g)}{_{ab}^{ba}}\(\xi\)\big).
\end{multline*}
By using these identities together with \eqref{Pr1Step3Eq5}--\eqref{Pr1Step3Eq7} and observing that
$$\(\nabla S\otimes g,\nabla^3U_\mu\)=-\Delta\(\nabla S,\nabla U_\mu\)-2\(\nabla^2S,\nabla^2U_\mu\)-\(\nabla^3S,\nabla U_\mu\otimes g\),$$
and
$$\tensor{(\nabla^2 S)}{_a^a}\(\xi\)=-\frac{1}{6}\left|\Weyl\(\xi\)\right|^2=-\frac{1}{6n}\left|\Weyl\(\xi\)\right|^2\tensor{g}{_a^a}\(\xi\),$$
we obtain that for $k\ge3$,
\begin{align}
&\int_M\(T_3,\nabla^2U_\mu\)\Delta^{k-3}U_\mu\,dv_g\nonumber\\
&\qquad=\frac{n+3k+1}{36n\(n-1\)}\left|\Weyl\(\xi\)\right|^2\int_M\(\Delta U_\mu\)\Delta^{k-3}U_\mu\,dv_g+\left\{\begin{aligned}&\bigO\(\mu^4\)&&\text{if }n=2k+4\\&\smallo\(\mu^4\)&&\text{if }n>2k+4\end{aligned}\right.\nonumber\\
&\qquad=\frac{n+3k+1}{36n\(n-1\)}\left|\Weyl\(\xi\)\right|^2\int_MU_\mu\Delta^{k-2}U_\mu\,dv_g+\left\{\begin{aligned}&\bigO\(\mu^4\)&&\text{if }n=2k+4\\&\smallo\(\mu^4\)&&\text{if }n>2k+4\end{aligned}\right.\label{Pr1Eq6}
\end{align}
and
\begin{align}
&\int_M\(T_4,\nabla^3U_\mu\)\Delta^{k-3}U_\mu\,dv_g=-\frac{k+1}{\(n-1\)\(n+2\)}\int_M\bigg(\Delta\(\nabla S,\nabla U_\mu\)\nonumber\\
&\quad+\frac{1}{3n}\left|\Weyl\(\xi\)\right|^2\Delta U_\mu+\(\nabla^3S,\nabla U_\mu\otimes g\)\bigg)\Delta^{k-3}U_\mu\,dv_g+\left\{\begin{aligned}&\bigO\(\mu^4\)&&\text{if }n=2k+4\\&\smallo\(\mu^4\)&&\text{if }n>2k+4\end{aligned}\right.\nonumber\allowdisplaybreaks\\
&=-\frac{k+1}{3n\(n-1\)\(n+2\)}\int_M\(\left|\Weyl\(\xi\)\right|^2U_\mu+3n\(\nabla S,\nabla U_\mu\)\)\Delta^{k-2}U_\mu\,dv_g\nonumber\\
&\quad-\frac{k+1}{\(n-1\)\(n+2\)}\int_MU_\mu\Delta^{k-3}\(\nabla^3S,\nabla U_\mu\otimes g\)dv_g+\left\{\begin{aligned}&\bigO\(\mu^4\)&&\text{if }n=2k+4\\&\smallo\(\mu^4\)&&\text{if }n>2k+4,\end{aligned}\right.\label{Pr1Eq7}
\end{align}
and for $k\ge4$,
\begin{align}
&\int_M\(T_5,\nabla^4U_\mu\)\Delta^{k-4}U_\mu\,dv_g\nonumber\\
&\quad=-\frac{\(k+1\)\left|\Weyl\(\xi\)\right|^2}{10n\(n-1\)\(n+2\)}\int_M\(\Delta^2U_\mu\)\Delta^{k-4}U_\mu\,dv_g+\left\{\begin{aligned}&\bigO\(\mu^4\)&&\text{if }n=2k+4\\&\smallo\(\mu^4\)&&\text{if }n>2k+4\end{aligned}\right.\nonumber\\
&\quad=-\frac{\(k+1\)\left|\Weyl\(\xi\)\right|^2}{10n\(n-1\)\(n+2\)}\int_MU_\mu\Delta^{k-2}U_\mu\,dv_g+\left\{\begin{aligned}&\bigO\(\mu^4\)&&\text{if }n=2k+4\\&\smallo\(\mu^4\)&&\text{if }n>2k+4.\end{aligned}\right.\label{Pr1Eq8}
\end{align}
It follows from \eqref{Pr1Eq4} and \eqref{Pr1Eq6}--\eqref{Pr1Eq8} that
\begin{align}\label{Pr1Eq9}
&\int_MU_\mu P_{2k}U_\mu\,dv_g=\int_MU_\mu\Delta^kU_\mu\,dv_g+k\int_M\bigg(\bigg(\(k-1\)J_2+\Delta J_1\nonumber\\
&+\(\frac{\(k-1\)\(k-2\)\(n+3k+1\)}{36n\(n-1\)}-\frac{\(k-1\)\(k-2\)\(k+1\)\(3k+1\)}{30n\(n-1\)\(n+2\)}\)\left|\Weyl\(\xi\)\right|^2\bigg)U_\mu\allowdisplaybreaks\nonumber\\
&+\(\(k-1\)T_1-2\nabla J_1-\frac{\(k+1\)\(k-1\)\(k-2\)}{\(n-1\)\(n+2\)}\nabla S,\nabla U_\mu\)\nonumber\\
&+\(\(k-1\)T_2-J_1g,\nabla^2U_\mu\)\bigg)\Delta^{k-2}U_\mu\,dv_g\allowdisplaybreaks\nonumber\\
&+\int_MU_\mu\(ZU_\mu-\frac{k\(k+1\)\(k-1\)\(k-2\)}{\(n-1\)\(n+2\)}\Delta^{k-3}\(\nabla^3S,\nabla U_\mu\otimes g\)\)\,dv_g\nonumber\\
&+\left\{\begin{aligned}&\bigO\(\mu^4\)&&\text{if }n=2k+4\\&\smallo\(\mu^4\)&&\text{if }n>2k+4.\end{aligned}\right.
\end{align}
Straightforward computations together with \eqref{Eq4}, \eqref{Eq5} and \eqref{Pr1Eq5} give
\begin{align*}
&J_1\(\xi\)=0\quad\text{and}\quad\Delta J_1\(\xi\)=\frac{n-2}{24\(n-1\)}\left|\Weyl\(\xi\)\right|^2,\allowdisplaybreaks\\
&J_2\(\xi\)=-\frac{3n+2k-4}{144\(n-1\)}\left|\Weyl\(\xi\)\right|^2,\allowdisplaybreaks\\
&\tensor{(T_1)}{^a_{;a}}\(\xi\)=-\frac{3n+4k-2}{72\(n-1\)}\left|\Weyl\(\xi\)\right|^2,\allowdisplaybreaks\\
&\tensor{(T_2)}{_a^a}\(\xi\)=0
\end{align*}
and
$$\tensor{(T_2)}{_a^a_{;b}^b}\(\xi\)=\tensor{(T_2)}{^{ab}_{;ab}}\(\xi\)=\tensor{(T_2)}{^{ab}_{;ba}}\(\xi\)=-\frac{k+1}{18\(n-1\)}\left|\Weyl\(\xi\)\right|^2.$$
By using these identities together with \eqref{Pr1Step2Eq1}, \eqref{Pr1Step3Eq2}--\eqref{Pr1Step3Eq4}, \eqref{Pr1Eq2} and \eqref{Pr1Eq9}, we obtain that \eqref{Pr1Eq1} holds true with $C\(n,k\)$ defined as
\begin{align}\label{Pr1Eq10}
&C\(n,k\):=\frac{\(n-3\)\(n-5\)!\,k!}{16\(k-1\)\(n-2k-2\)}\nonumber\\
&\times\sum_{l=k-2}^{2k-4}\frac{l!}{\(l-k+2\)!\(2k-l-4\)!\(n+l-2k+1\)!}\bigg(8\(n+l-2k\)\(n+l-2k+1\)\allowdisplaybreaks\nonumber\\
&\quad\times\bigg(\frac{\(k-1\)\(3n+2k-4\)}{144}-\frac{n-2}{24}-\frac{\(k-1\)\(k-2\)\(n+3k+1\)}{36n}\nonumber\\
&\quad+\frac{\(k-1\)\(k-2\)\(k+1\)\(3k+1\)}{30n\(n+2\)}\bigg)+4\(n-2k\)\(n+l-2k+1\)\bigg(\frac{n-2}{12}\nonumber\\
&\quad-\frac{\(k-1\)\(3n+4k-2\)}{72}+\frac{\(k+1\)\(k-1\)\(k-2\)}{6\(n+2\)}\bigg)\allowdisplaybreaks\nonumber\\
&\quad+\(n-2k\)\bigg(\frac{\(k+1\)\(k-1\)\(n-2k-2l+4\)}{18}\nonumber\\
&\quad-\frac{\(n-2\)\(2\(n-2k+2\)-n\(n+2l-2k\)\)}{24}\bigg)\bigg)\allowdisplaybreaks\nonumber\\
&\quad\times\Beta\(\frac{n}{2}-k-1,l+1\)^{-1}\left\{\begin{aligned}&2\chi_{\left\{l=k-2\right\}}&&\text{if }n=2k+4\\&\Beta\(\frac{n}{2}+l-2k,2k-l-2\)&&\text{otherwise}\end{aligned}\right.\allowdisplaybreaks\nonumber\\
&=\frac{\(n-3\)\(n-5\)!\,k!}{5760n\(n+2\)\(k-1\)\(n-2k-2\)}\nonumber\\
&\quad\times\sum_{l=k-2}^{2k-4}\frac{l!\,c\(n,k,l\)}{\(l-k+2\)!\(2k-l-4\)!\(n+l-2k+1\)!}\nonumber\\
&\quad\times\Beta\(\frac{n}{2}-k-1,l+1\)^{-1}\left\{\begin{aligned}&2\chi_{\left\{l=k-2\right\}}&&\text{if }n=2k+4\\&\Beta\(\frac{n}{2}+l-2k,2k-l-2\)&&\text{otherwise,}\end{aligned}\right.
\end{align}
where
\begin{align*}
&c\(n,k,l\):=4\(n+l-2k\)\(n+l-2k+1\)(5n\(n+2\)\(k-1\)\(3n+2k-4\)\\
&\quad-30n\(n+2\)\(n-2\)-20\(n+2\)\(k-1\)\(k-2\)\(n+3k+1\)\allowdisplaybreaks\\
&\quad+24\(k-1\)\(k-2\)\(k+1\)\(3k+1\))+20n\(n-2k\)\(n+l-2k+1\)\allowdisplaybreaks\\
&\quad\times\(6\(n+2\)\(n-2\)-\(n+2\)\(k-1\)\(3n+4k-2\)+12\(k+1\)\(k-1\)\(k-2\)\)\allowdisplaybreaks\\
&\quad+5n\(n+2\)\(n-2k\)(4\(k+1\)\(k-1\)\(n-2k-2l+4\)\\
&\quad-3\(n-2\)\(2\(n-2k+2\)-n\(n+2l-2k\)\)).
\end{align*}
By letting $k:=3+a$, $n:=2k+4+b$ and $l:=k-2+c$ and using the software {\it Maple} to expand the expression of $c\(n,k,l\)$, we then obtain
\begin{align}\label{Pr1Eq11}
&c\(n,k,l\)=4(15ab^3+1200a^2b+3880ab+1920+10656a+480b+4528a^2+624a^3\nonumber\\
&\quad+40b^2+450ab^2+80a^3b+80a^2b^2+32a^4)c^2+2(71552a^2b+414912a+500a^2b^3\allowdisplaybreaks\nonumber\\
&\quad+247984ab+31840a^3+53660ab^2+3200a^4+640a^3b^2+11020b^3+150ab^4\allowdisplaybreaks\nonumber\\
&\quad+128a^5+9056a^3b+660b^4+161440a^2+448a^4b+15b^5+311040b+10520a^2b^2\allowdisplaybreaks\nonumber\\
&\quad+4830ab^3+426240+85840b^2)c+128a^6+576a^5b+1088a^4b^2+1020a^3b^3\allowdisplaybreaks\nonumber\\
&\quad+560a^2b^4+150ab^5+15b^6+3904a^5+14720a^4b+21896a^3b^2+15940a^2b^3\allowdisplaybreaks\nonumber\\
&\quad+5640ab^4+720b^5+49408a^4+149280a^3b+167032a^2b^2+81120ab^3+13780b^4\allowdisplaybreaks\nonumber\\
&\quad+332096a^3+754720a^2b+563824ab^2+134240b^3+1250304a^2+1900224ab\nonumber\\
&\quad+704640b^2+2499840a+1900800b+2073600.
\end{align}
Since all the coefficients in this expression are positive, it follows that $C\(n,k\)$ is positive whenever $k\ge3$, $n\ge2k+4$ and $l\ge k-2$. Furthermore, in the case where $k=2$ and $l=0$, we find
$$c\(n,2,0\)=5n(n+2)(n-4)^2(n^2-4n-4)>0\quad\forall n\ge8.$$
Therefore, in all cases, we find that $C\(n,k\)$ is positive. This ends the proof of Proposition~\ref{Pr1}.
\endproof

We can now prove Theorem~\ref{Th1} by using Proposition~\ref{Pr1}.

\proof[Proof of Theorem~\ref{Th1} in the case where $n\ge2k+4$]
Let $\xi\in M$ be a maximal point of $f$ and $\widetilde{g}=\varphi^{4/\(n-2\)}g$ be a conformal metric to $g$ such that $\varphi\(\xi\)=1$ and $\det\widetilde{g}\(x\)=1$ for all $x$ in a neighborhood of the point $\xi$. Notice that since $\xi$ is a maximal point of $f$, if $\Delta_gf\(\xi\)=0$, then $\nabla^jf=0$ for all $j\in\left\{1,2,3\right\}$. In particular, since $\varphi\(\xi\)=1$, it follows that
\begin{equation}\label{Th1Eq4}
\Delta_{\widetilde{g}}f\(\xi\)=0\quad\text{and}\quad\Delta_{\widetilde{g}}^2f\(\xi\)=\Delta_g^2f\(\xi\),
\end{equation}
where $\Delta_g$ and $\Delta_{\widetilde{g}}$ are the Laplace--Beltrami operators with respect to the metrics $g$ and $\widetilde{g}$, respectively, and the covariant derivatives, the Ricci tensor and the multiple inner product in the right-hand side of the second identity are with respect to the metric $g$. Let $c\(n,k\)$ be the constant defined as 
\begin{equation}\label{Th1Eq5}
c\(n,k\):=\left\{\begin{aligned}&0&&\text{if }n=2k+4\\&\frac{\(n-2k\)\(2k-1\)!}{8n\(n-2\)\(n-4\)C\(n,k\)}\Beta\(\frac{n}{2}-k,2k\)^{-1}&&\text{if }n>2k+4,\end{aligned}\right.
\end{equation}
where $C\(n,k\)$ is as in \eqref{Pr1Eq1} (see also \eqref{Pr1Eq10}). By applying Proposition~\ref{Pr1} together with \eqref{Th1Eq4} and the fact that $\left|W\right|$ is conformally invariant, we then obtain that if \eqref{Th1Eq1} and \eqref{Th1Eq2} hold true, then
\begin{equation}\label{Th1Eq6}
\inf_{u\in C^{2k}\(M\)\backslash\left\{0\right\}}\hspace{-2pt}I_{k,f,\widetilde{g}}\(u\)<\omega_n^{\frac{2k}{n}}\(2k-1\)!\Beta\(\frac{n}{2}-k,2k\)^{-1}\hspace{-2pt}\Big(\max_{x\in M} f\(x\)\Big)^{-\frac{n-2k}{n}}.
\end{equation}
On the other hand, by conformal invariance of the operator $P_{2k}$, we obtain
\begin{equation}\label{Th1Eq7}
\inf_{u\in C^{2k}\(M\)\backslash\left\{0\right\}}I_{k,f,\widetilde{g}}\(u\)=\inf_{u\in C^{2k}\(M\)\backslash\left\{0\right\}}I_{k,f,g}\(u\).
\end{equation}
By putting together \eqref{Th1Eq6} and \eqref{Th1Eq7} and applying Theorem~3 of Mazumdar~\cite{Maz}, we then obtain that the conclusions of Theorem~\ref{Th1} hold true. 
\endproof

\section{The remaining cases}\label{Sec3}

This section is devoted to the proof of Theorem~\ref{Th1} in the remaining case where $2k+1\le n\le 2k+3$ together with the following result in the case where $g$ is conformally flat in some open subset of the manifold:

\begin{theorem}\label{Th3}
Let $k\ge1$ be an integer, $\(M,g\)$ be a smooth, closed Riemannian manifold of dimension $n\ge 2k+1$ and $f$ be a smooth positive function in $M$. Assume that $Y_{2k}>0$ and there exists a maximal point $\xi$ of $f$ such that $m\(\xi\)>0$ (see \eqref{Sec3Eq2} for the definition of the mass), $\nabla^jf\(\xi\)=0$ for all $j\in\left\{1,\dotsc,n-2k\right\}$ and $g$ is conformally flat in some neighborhood of the point $\xi$. Then there exists a nontrivial solution $u\in C^{2k}\(M\)$ to the equation \eqref{Eq1}, which minimizes the energy functional \eqref{Eq2}. If moreover the Green's function of the operator $P_{2k}$ is positive, then $u$ is positive, which implies that the $Q$-curvature of order $2k$ of the metric $u^{4/\(n-2k\)}g$ is equal to $\frac{2}{n-2k}f$. 
\end{theorem}

Notice that Theorem~\ref{Th2} is now a direct consequence of Theorems~\ref{Th1} and~\ref{Th3}.

\smallskip
Throughout this section, we fix a point $\xi\in M$ and assume that $2k+1\le n\le 2k+3$ or $g$ is conformally flat in some neighborhood of $\xi$. In these cases, our proofs are based on the  method of Schoen \cite{Sch1} for the resolution of the remaining cases of the Yamabe problem, which has been extended to the $k=2$ case by Gursky and Malchiodi~\cites{GurMal} and Hang and Yang~\cites{HangYang1,HangYang2}. We  consider a family of global test functions involving the Green's function and derive an expression for  the energy functional $I_{k,f,g}$ (see \eqref{Eq2}) associated with the equation \eqref{Eq1}. Then, analogously as in the case $n\geq 2k+4$, by using the expansion obtained in Proposition~\ref{Pr1}, we obtain the existence of a nontrivial solution to the equation \eqref{Eq1} under a positivity assumption on the mass of the operator $P_{2k}$.

\smallskip
We now discuss the definition of the mass. By applying a conformal change of metric, we may assume that
\begin{equation}\label{Sec3Eq1}
\left\{\begin{aligned}&\text{$g$ satisfies \eqref{Eq3} in some neighborhood $\Omega$ of $\xi$ if $2k+1\le n\le2k+3$}\\&\text{$g$ is flat in some neighborhood $\Omega$ of $\xi$ if $n\ge2k+4$.}\end{aligned}\right.
\end{equation}
Then, in the geodesic normal coordinates at $\xi$ determined by $g$, the Green's function $G_{2k}\(x\):=G_{2k}\(x,\xi\)$ of the operator $P_{2k}$ has the expansion 
\begin{equation}\label{Sec3Eq2}
G_{2k}\(x\)=b_{n,k}\,d_g\(x,\xi\)^{2k-n}+m\(\xi\)+\smallo\(1\)
\end{equation}
as $x\to\xi$ (see Lee and Parker~\cite{LeePar} for $k=1$ and Michel~\cite{Mic} for $k\ge2$), where $m\(\xi\)\in M$ is called the mass of the operator $P_{2k}$ at the point $\xi$ and the constant $b_{n,k}$ is defined as 
$$b_{n,k}^{-1}:=2^{k-1}\(k-1\)!\(n-2\)\(n-4\)\dotsm\(n-2k\)\omega_{n-1}.$$
It is important to point out that the sign of $m\(\xi\)$ does not depend on our choice of conformal metric (see Michel~\cite{Mic}*{Th\'eor\`eme~3.1}). 

\smallskip
Now that the mass is defined, we consider the regular part of the Green's function, which plays a crucial role in the proofs of our theorems. It follows from \eqref{Sec3Eq2} that there exists a continuous function $h_{2k}$ in $M$ such that $h_{2k}\(\xi\)=m\(\xi\)$ and
\begin{equation}\label{Sec3Eq3}
G_{2k}\(x\)=b_{n,k}\,d_{g}\(x,\xi\)^{2k-n}+h_{2k}\(x\)\qquad\forall x\in M\backslash\left\{\xi\right\}.
\end{equation}
Furthermore, we have that $h_{2k}\in C^\infty\(\Omega\)$ in the case where $g$ is flat in $\Omega$ and $h_{2k}\in W^{2k,p}\(\Omega\)$ for all $p\in\[1,n/\(n-4\)\)$ if $n\ge5$ and $p\in\[1,\infty\)$ if $n\in\left\{3,4\right\}$ in the case where $2k+1\le n\le2k+3$ and $g$ satisfies \eqref{Eq3} in $\Omega$. This follows from classical elliptic regularity theory (see Agmon, Douglis and Nirenberg~\cite{AgmDouNir}) together with the fact that $P_{2k}h_{2k}=\Delta^kh_{2k}=0$ in $\Omega$ in the case where $g$ is flat in $\Omega$ and $P_{2k}h_{2k}=\bigO(d_g\(\cdot,\xi\)^{4-n})$ in $\Omega$ in the case where $g$ satisfies \eqref{Eq3} in $\Omega$ (see~\cite{Mic}*{Lemme~2.2}).
 
\smallskip
For every $\mu>0$, letting $\chi$ and $U_\mu$ be as in Section~\ref{Sec2}, we consider the test functions $V_\mu$ defined as
\begin{equation}\label{Sec3Eq4}
V_\mu\(x\):= U_\mu\(x\) + b^{-1}_{n,k}\,\mu^{\frac{n-2k}{2}} \(\chi(d_g\(x,\xi\))h_{2k}\(x\)+\(1-\chi\(d_g\(x,\xi\)\)\)G_{2k}\(x\)\)\hspace{-1pt}
\end{equation}
for all $x\in M$. Note that $V_\mu\in W^{2k,2n/\(n+2k\)}\(M\)$ so that in particular the integral $\int_{M}V_\mu P_{2k}V_\mu\,dv_g$ is well defined. We then obtain the following:

\begin{proposition}\label{Pr2}
Let $k\ge1$ be an integer, $\(M,g\)$ be a smooth, closed Riemannian manifold of dimension $n\ge2k+1$, $f$ be a smooth positive function in $M$ and $\xi$ be a point in $M$ such that $\nabla^jf\(\xi\)=0$ for all $j\in\left\{1,\dotsc,n-2k\right\}$ and \eqref{Sec3Eq1} holds true. Let $I_{k,f,g}$ be as in \eqref{Eq2} and $V_\mu$ be as in \eqref{Sec3Eq4}. Then, as $\mu\to0$,
\begin{multline}\label{Pr2Eq}
I_{k,f,g}\(V_\mu\)=\omega_n^{\frac{2k}{n}}\(2k-1\)!\Beta\(\frac{n}{2}-k,2k\)^{-1}f\(\xi\)^{-\frac{n-2k}{n}}\\
\times\Bigg(1-b_{n,k}^{-1}\Beta\(\frac{n}{2},\frac{n}{2}\)^{-1}\Beta\(\frac{n}{2},k\)m\(\xi\) \mu^{n-2k}+\smallo\(\mu^{n-2k}\)\Bigg).
\end{multline}
\end{proposition}

\proof[Proof of Proposition~\ref{Pr2}] The first step in the proof is as follows:

\begin{step}\label{Pr2Step1}
Assume that $g$ satisfies \eqref{Sec3Eq1} for some point $\xi\in M$. Then, as $\mu \to 0$,
\begin{multline}\label{Pr2Step1Eq1}
\int_{M}V_\mu P_{2k}V_\mu\,dv_g= 2^{2k-n}\omega_{n}(2k-1)!\Beta\(\frac{n}{2}-k,2k\)^{-1}\\
\times\(1+b^{-1}_{n,k}\Beta\(\frac{n}{2},\frac{n}{2}\)^{-1}\Beta\(\frac{n}{2},k\)m\(\xi\)\mu^{n-2k}+\smallo\(\mu^{n-2k}\)\).
\end{multline}
\end{step}

\proof[Proof of Step~\ref{Pr2Step1}]
We write
\begin{equation*}
V_\mu\(x\)=b_{n,k}^{-1}\,\mu^{\frac{n-2k}{2}}G_{2k}\(x\)+\Big(\underbrace{  U_\mu\(x\) - \mu^{\frac{n-2k}{2}}  \chi\(d_g\(x,\xi\)\)d_g\(x,\xi\)^{2k-n} }_{W_\mu\(x\)} \Big)
\end{equation*}
for all $x\in M\backslash\left\{\xi\right\}$. Straightforward estimates give
\begin{align}\label{Pr2Step1Eq2}
&\int_{M}V_\mu P_{2k}V_\mu\,dv_g=\int_{B\(\xi,2r_0\)} V_\mu P_{2k}W_\mu\,dv_g\nonumber\\
&\quad=\int_{B\(\xi,r_0\)} \(U_{\mu}+ b_{n,k}^{-1}\,\mu^{\frac{n-2k}{2}} h_{2k} \)P_{2k}W_{\mu}\,dv_{g}\nonumber\\
&\qquad+\bigO \Bigg( \mu^{\frac{n-2k}{2}} \int_{B\(\xi,2r_0\) \setminus B\(\xi,r_0\)}\left|P_{2k}W_{\mu}\right|\,dv_{g} \Bigg)\allowdisplaybreaks\nonumber\\
&\quad=\int_{B\(\xi,r_0\)} U_{\mu}P_{2k}W_{\mu}\,dv_{g}+ b_{n,k}^{-1}\,\mu^{\frac{n-2k}{2}} \int_{B\(\xi,r_0\)}  h_{2k}P_{2k}W_{\mu}\,dv_{g}\nonumber\\
&\qquad+\bigO \Bigg( \mu^{n-2k} \sum \limits_{\left|\alpha\right| \leq 2k} \int_{B\(0,2r_0\) \setminus B\(0,r_0\)} \Big|\partial^{\alpha}\Big[\big(\mu^{2}+\left|x\right|^{2} \big)^{\frac{2k-n}{2}}-\left|x\right|^{2k-n}\Big]\Big|\,dx \Bigg)\nonumber\\
&\quad=\int_{B\(\xi,r_0\)} U_{\mu}P_{2k}W_{\mu}\,dv_{g}+b_{n,k}^{-1}\,\mu^{\frac{n-2k}{2}} \int_{B\(\xi,r_0\)}  h_{2k}P_{2k}W_{\mu}\,dv_{g}+\smallo(\mu^{n-2k}).
\end{align}
We claim that
\begin{equation}\label{Pr2Step1Eq3}
\left| P_{2k}W_{\mu}-\Delta^{k}U_{\mu}\right|\leq C \mu^{\frac{n-2(k-2)}{2}}d_g\(x,\xi\)^{2-n} \quad\forall x\in B\(\xi,r_{0}\)\setminus\left\{\xi\right\}
\end{equation}
for some constant $C$ independent of $x$, $\mu$ and $\xi$. Assuming \eqref{Pr2Step1Eq3} and proceeding as in Step~\ref{Pr1Step3}, it then follows from \eqref{Pr2Step1Eq2} that
\begin{align*}
&\int_{M}V_\mu P_{2k}V_\mu\,dv_g=\int_{B\(\xi,r_0\)}U_{\mu}\Delta^{k}U_{\mu}\,dv_{g}+b_{n,k}^{-1}\,\mu^{\frac{n-2k}{2}} \int_{B\(\xi,r_0\)}  h_{2k}\Delta^{k}U_{\mu}\,dv_{g}\\
&\qquad +\bigO\(\mu^{n-2(k-1)}\int_{B\(0,r_{0}\)}\left|x\right|^{2-n}\big(\mu^{2}+\left|x\right|^{2} \big)^{\frac{2k-n}{2}}\,dx \) +\smallo\(\mu^{n-2k}\)\allowdisplaybreaks\\
&\quad =\int_{B\(0,r_0\)}  \widetilde{U}_{\mu}\Delta^{k}_{0}  \widetilde{U}_{\mu}\,dx+ b_{n,k}^{-1}\,\mu^{\frac{n-2k}{2}}\int_{B\(0,r_0\)}  h_{2k}\(\exp_\xi x\)\Delta^{k}_{0} \widetilde{U}_{\mu}\,dx +\smallo\(\mu^{n-2k}\)\allowdisplaybreaks\\
&\quad = 2^{2k-1}\big( 2k-1\big)!\,\omega_{n-1}  \Beta\(\frac{n}{2}-k,2k\)^{-1} \Bigg(\int_0^{\(r_0/\mu\)^2}\frac{r^{\frac{n-2}{2}}}{\(1+r\)^{n}}\,dr\\
&\qquad+b_{n,k}^{-1}m\(\xi\)\mu^{n-2k}\int_0^{\(r_0/\mu\)^2}\frac{r^{\frac{n-2}{2}}}{\(1+r\)^{\frac{n+2k}{2}}}dr\Bigg) +\smallo\(\mu^{n-2k}\)\allowdisplaybreaks\\
&\quad= 2^{2k-n}\(2k-1\)!\,\omega_nB\(\frac{n}{2}-k,2k\)^{-1}\\
&\qquad\times\(1+b^{-1}_{n,k}\Beta\(\frac{n}{2},\frac{n}{2}\)^{-1}\Beta\(\frac{n}{2},k\)m\(\xi\)\mu^{n-2k}+\smallo\(\mu^{n-2k}\)\).
\end{align*}
Therefore, it remains to prove \eqref{Pr2Step1Eq3} to complete the proof of Step~\ref{Pr2Step2}. Notice that \eqref{Pr2Step1Eq3} is clearly satisfied with $C=0$ in the case where $n\ge 2k+4$ and $g$ is flat in $\Omega$. Therefore, we may assume in what follows that we are in the case where $2k+1\le n\le 2k+3$ and $g$ satisfies \eqref{Eq3} in $\Omega$. By using \eqref{Pr1Step1Eq1}, we obtain 
\begin{multline}\label{Pr2Step1Eq4}
P_{2k}= \Delta^k+ k\Delta^{k-1}\(J_1\cdot\) + k\(k-1\)\Delta^{k-2}\(\(T_1,\nabla \)+\(T_2,\nabla^{2}\)\)\\+k\(k-1\)\(k-2\)\Delta^{k-3}\(T_4,\nabla^3\)+Z,
\end{multline}
where $Z$ is a smooth linear operator of order less than $2k-3$ if $k\ge2$, $Z:=0$ if $k=1$. By induction, one can check that 
\begin{align}
&\Delta^{k-1}\(J_1\cdot\)=J_1\Delta^{k-1}-2\(k-1\)\(\nabla J_1,\nabla\Delta^{k-2}\)+\smallo^{2k-3},\label{Pr2Step1Eq5}\allowdisplaybreaks\\
&\Delta^{k-2}\(T_1,\nabla\)=\(T_1,\nabla\Delta^{k-2}\)+\smallo^{2k-3},\label{Pr2Step1Eq6}\allowdisplaybreaks\\
&\Delta^{k-2}\(T_2,\nabla^2\)=\(T_2,\nabla^2\Delta^{k-2}\)-2\(k-2\)\(\nabla T_2,\nabla^3\Delta^{k-3}\)+\smallo^{2k-3}\label{Pr2Step1Eq7}
\end{align}
and
\begin{equation}\label{Pr2Step1Eq8}
\Delta^{k-3}\(T_4,\nabla^3\)=\(T_4,\nabla^{3}\Delta^{k-3}\)+\smallo^{2k-3},
\end{equation}
where $\smallo^{2k-3}$ is as in the proof of Step~\ref{Pr1Step1}. It follows from \eqref{Eq4}, \eqref{Eq5} and \eqref{Pr2Step1Eq4}--\eqref{Pr2Step1Eq8} that
\begin{align}\label{Pr2Step1Eq9}
P_{2k}W_{\mu}&=\Delta^k W_{\mu}+\frac{2k\(k-1\)\(k+1\)}{3\(n-2\)}\big(\(\Ricci,\nabla^2\Delta^{k-2}W_\mu\)-\(\delta\Ricci,\nabla\Delta^{k-2}W_\mu\)\nonumber\\
&\quad-\(k-2\)\(\nabla\Ricci,\nabla^3\Delta^{k-2}W_\mu\)\big)+\bigO\Big(d_g\(\cdot,\xi\)^2\left|\nabla^{2k-2}W_\mu\right|\nonumber\\
&\quad+d_g\(\cdot,\xi\)\left|\nabla^{2k-3}W_\mu\right|+\sum_{j=0}^{2k-4}\left|\nabla^jW_\mu\right|\Big)
\end{align}
in $M\backslash\left\{\xi\right\}$, uniformly with respect to $\mu$ and $\xi$. By using geodesic normal coordinates together with \eqref{Eq4} and a Taylor expansion, it follows from \eqref{Pr2Step1Eq9} that
\begin{multline}\label{Pr2Step1Eq10}
\(P_{2k}W_{\mu}-\Delta^k U_{\mu}\)\(\exp_\xi x\)=\(P_{2k}W_{\mu}-\Delta^k W_{\mu}\)\(\exp_\xi x\)\\
=\bigO\Big(\left|x\right|^2\big|\nabla^{2k-2}\widetilde{W}_\mu\(x\)\big|+\left|x\right|\big|\nabla^{2k-3}\widetilde{W}_\mu\(x\)\big|+\sum_{j=0}^{2k-4}\big|\nabla^j\widetilde{W}_\mu\(x\)\big|\Big)
\end{multline}
uniformly with respect to $x\in B\(0,r_0\)\backslash\left\{0\right\}$, $\mu$ and $\xi$, where
$$\widetilde{W}_\mu\(x\):=\mu^{\frac{2k-n}{2}}U\(x/\mu\)-\mu^{\frac{n-2k}{2}}\left|x\right|^{2k-n}.$$
Similarly as in \eqref{Pr1Step3Eq15}, we obtain that for every $j\in\N$, there exists a constant $C_j$ independing of $x$ and $\mu$ such that 
\begin{align}\label{Pr2Step1Eq11}
\big|\nabla^j\widetilde{W}_\mu\(x\)\big|&=\mu^{\frac{2k-n-2j}{2}}\left|\nabla^jW\(x/\mu\)\right|\nonumber\\
&\le C_j\sum_{m=0}^{\lfloor j/2\rfloor}\mu^{\frac{2k-n-4j+4m}{2}}r^{\frac{j-2m}{2}}\left|\partial^{j-m}_rW\(r/\mu^2\)\right|
\end{align}
for all $x\in B\(0,r_0\)\backslash\left\{0\right\}$, where $r:=\left|x\right|^2$ and 
$$W\(x\)=W\(r\):=\(1+r\)^{\(2k-n\)/2}-r^{\(2k-n\)/2}.$$
Furthermore, it is easy to see that 
\begin{equation}\label{Pr2Step1Eq12}
\left|\partial^{j}_rW\(r\)\right|\le C_j'r^{\frac{2k-n-2j-2}{2}}
\end{equation}
for some constant $C'_j$ independent of $r$. We then obtain \eqref{Pr2Step1Eq3} by putting together \eqref{Pr2Step1Eq10}--\eqref{Pr2Step1Eq12}. This completes the proof of Step~\ref{Pr2Step1}.
\endproof

\begin{step}\label{Pr2Step2}
Assume that $g$ satisfies \eqref{Sec3Eq1} for some point $\xi\in M$. Let $f$ be a smooth function $f$ in $M$ such that $\nabla^jf\(\xi\)=0$ for all $j\in\left\{1,\dotsc,n-2k\right\}$. Then,  as $\mu\to0$,
\begin{multline}\label{Pr2Step2Eq1}
\int_{M}f\left|V_\mu\right|^{2^*_k}dv_g=\frac{\omega_n}{2^n}\bigg(f\(\xi\)+2^*_k\,b^{-1}_{n,k}\Beta\(\frac{n}{2},\frac{n}{2}\)^{-1}\Beta\(\frac{n}{2},k\)f\(\xi\)m\(\xi\)\mu^{n-2k}\\
+\smallo\(\mu^{n-2k}\)\bigg).
\end{multline}
\end{step}

\proof[Proof of Step~\ref{Pr2Step2}]
By using a Taylor expansion together with straightforward estimates, we obtain
\begin{align}
&\int_{M}f \left|V_\mu\right|^{2^*_k}dv_g= \int_{B\(\xi, r_{0}\)}f \big| U_{\mu}+ b_{n,k}^{-1}\,\mu^{\frac{n-2k}{2}}h_{2k}\big|^{2^*_k}dv_g + \bigO\(\mu^{n}\)\notag\\
&=f\(\xi\)\int_{B\(\xi, r_{0}\)} U_{\mu}^{2^*_k}dv_g + 2^*_k\,b_{n,k}^{-1}\,\mu^{\frac{n-2k}{2}} \int_{B\(\xi, r_{0}\)}fh_{2k}U_{\mu}^{2^*_k-1}dv_g\notag\\
&\quad+ \bigO\(  \int_{B\(0, r_{0}\)}\Bigg(\left|x\right|^{n-2k+1}\( \frac{\mu}{\mu^{2}+\left|x\right|^{2}}\)^{n}+ \mu^{n-2k}\( \frac{\mu}{\mu^{2}+\left|x\right|^{2}}\)^{2k}\Bigg)dx+\mu^n\)\notag\allowdisplaybreaks\\
&=\frac{\omega_n}{2^n}f\(\xi\) + \frac{n\,\omega_{n-1}}{n-2k}\,b_{n,k}^{-1}f\(\xi\) m\(\xi\) \mu^{n-2k} \int_0^{\(r_0/\mu\)^2}\hspace{-10pt}\frac{r^{\frac{n-2}{2}}dr}{\(1+r\)^{\frac{n+2k}{2}}}+  \smallo\(\mu^{n-2k}\).\label{Pr2Step2Eq2}
\end{align}
Then \eqref{Pr2Step2Eq1} follows from \eqref{Pr1Step3Eq20},  \eqref{Pr1Step3Eq22} and \eqref{Pr2Step2Eq2}.
\endproof

We can now end the proofs of Proposition~\ref{Pr2} and Theorems~\ref{Th1} and~\ref{Th2}.

\proof[End of proof of Proposition~\ref{Pr2}]
We obtain \eqref{Pr2Eq} by putting together \eqref{Pr2Step1Eq1} and \eqref{Pr2Step2Eq1}.
\endproof

\proof[Proofs of Theorem~\ref{Th1} in the case where $2k+1\le n\le 2k+3$ and of Theorem~\ref{Th3}] 
Let $\xi\in M$ be a maximal point of $f$ such that $\nabla^jf\(\xi\)=0$ for all $j\in\left\{1,\dotsc,n-2k\right\}$ (notice that for a maximal point, this is equivalent to \eqref{Th1Eq1} in the case where $2k+1\le n\le 2k+3$). By applying Proposition~\ref{Pr2} together with a conformal change of metric, we then obtain that if $m\(\xi\)>0$, then there exists a function $V\in W^{2k,2n/\(n+2k\)}\(M\)\backslash\left\{0\right\}$ such that 
$$I_{k,f,g}\(V\)<\omega_n^{\frac{2k}{n}}\(2k-1\)!\Beta\(\frac{n}{2}-k,2k\)^{-1}\Big(\max_{x\in M} f\(x\)\Big)^{-\frac{n-2k}{n}}.$$ 
Notice that $W^{2k,2n/\(n+2k\)}\(M\)\hookrightarrow L^{2^*_k}\(M\)$ so that a density argument gives
$$\inf_{u\in C^{2k}\(M\)\backslash\left\{0\right\}}I_{k,f,g}\(u\)<\omega_n^{\frac{2k}{n}}\(2k-1\)!\Beta\(\frac{n}{2}-k,2k\)^{-1}\Big(\max_{x\in M} f\(x\)\Big)^{-\frac{n-2k}{n}}.$$ 
We can then conclude the proofs of Theorems~\ref{Th1} and~\ref{Th3} by applying Theorem~3 of Mazumdar~\cite{Maz}.
\endproof

\end{document}